\documentclass{article}

\usepackage{enumitem}
\usepackage{mathrsfs}
\usepackage{fancyhdr}
\usepackage[T1]{fontenc}
\usepackage{amssymb,amsmath,amsfonts,amsxtra,amsthm}
\usepackage{graphicx}
\usepackage{placeins}

\begin{document}

\newtheorem{Theo}{Theorem}
\newtheorem{Prop}{Proposition}
\newtheorem{Coro}{Corollary}
\newtheorem{Lem}{Lemma}

\newcommand{\NN}{\mathbb{N}}
\newcommand{\RR}{\mathbb{R}}
\newcommand{\PP}{\mathbb{P}}
\newcommand{\EE}{\mathbb{E}}
\newcommand{\ind}{ {{\rm 1}\hskip-2.2pt{\rm l}}}
\newcommand{\union}{\cup}
\newcommand{\inter}{\cap}
\newcommand{\vide}{\varnothing}
\newcommand{\Var}{\operatorname{Var}}
\newcommand{\Cov}{\operatorname{Cov}}
\newcommand{\CQFD}
{%
\mbox{}%
\nolinebreak%
\hfill%
\rule{2mm}{2mm}%
\medbreak%
\par%
}

\title{Uniform strong consistency of a frontier estimator using kernel regression on high order moments}
\author{St\'ephane Girard$^{(1)}$, Armelle Guillou$^{(2)}$ \& Gilles Stupfler$^{(3)}$}
\date{\small $^{(1)}$ Team Mistis, INRIA Rh\^one-Alpes \& LJK, Inovall\'ee, 655, av. de l'Europe, \\ Montbonnot, 38334 Saint-Ismier cedex, France \\
$^{(2)}$ Universit\'e de Strasbourg \& CNRS, IRMA, UMR 7501, 7 rue Ren\'e Descartes,\\ 67084 Strasbourg cedex, France \\
$^{(3)}$ Aix Marseille Universit\'e, CERGAM, EA 4225, 15-19 all\'ee Claude Forbin, \\ 13628 Aix-en-Provence Cedex 1, France}

\maketitle

{\bf Abstract.} We consider the high order moments estimator of the frontier of a random pair, introduced by Girard, S., Guillou, A., Stupfler, G. (2013). {\it Frontier estimation with kernel regression on high order moments}. In the present paper, we show that this estimator is strongly uniformly consistent on compact sets and its rate of convergence is given when the conditional cumulative distribution function belongs to the Hall class of distribution functions. \\



\section{Introduction} 
\noindent
Let $(X_1, \, Y_1), \ldots, \, (X_n, \, Y_n)$ be $n$ independent copies of a random pair $(X, \, Y)$ such that their common distribution has a support defined by $S=\left\{  (x, \, y) \in E \times\mathbb{R}; \ 0\leq y\leq g(x) \right\}$, where $E$ is a closed subset of $\RR^d$ having nonempty interior. The unknown function $g$ is called the frontier. In Girard {\it et al.} (2013), a new estimator of $g$ is introduced, based upon kernel regression on high order moments of the data:
\begin{equation}
\label{defest}
\frac{1}{\widehat{g}_n(x)} 
= \frac{1}{a p_n} \left[ ((a+1)p_n+1) \frac{\widehat{\mu}_{(a+1)p_n}(x)}{\widehat{\mu}_{(a+1)p_n+1}(x)} - (p_n+1) \frac{\widehat{\mu}_{p_n}(x)}{\widehat{\mu}_{p_n+1}(x)} \right]
\end{equation}
where $(p_n)$ is a nonrandom positive sequence such that $p_n \to \infty$, $a>0$ and 
$$
\widehat{\mu}_{p_n}(x)=\frac{1}{n} \sum_{i=1}^n Y_i^{p_n} \, K_{h_n}(x-X_i)
$$
is a kernel estimator of the conditional moment $m_{p_n}(x)=\EE(Y^{p_n} \, | \, X=x)$. Classically, $K$ is a probability density function on $\RR^d$, $K_h(u)=h^{-d} \, K(u/h)$ and $(h_n)$ is a nonrandom positive sequence such that $h_n \to 0$. From a practical point of view, the use of a small window-width $h_n$ allows to select the pairs $(X_i, \, Y_i)$ such that $X_i$ is close to $x$ while the use of the high power $p_n$ gives more weight to the $Y_i$ close to $g(x)$. Using high order moments was first suggested by Girard and Jacob (2008) in the case when $Y$ given $X$ is uniformly distributed. This approach was also used in Girard and Jacob (2009) to develop a local polynomial estimator. In Girard {\it et al.} (2013), the estimator~(\ref{defest}) was shown to be pointwise consistent and asymptotically normal. Our focus in the present paper is to examine its almost sure uniform properties.

\vskip1ex
\noindent
Uniform consistency results in frontier estimation are seldom available in the literature: we refer the reader to Geffroy (1964) for the uniform consistency of the blockwise maxima estimator when the conditional distribution function of $Y$ given $X$ is uniform and to Jacob and Suquet (1995) for the uniform consistency of a projection estimator when the observations are realizations of a Poisson process whose intensity is known. Neither of these papers provides the rate of uniform convergence of the estimator it studies. In the field of econometrics, where the frontier function is assumed to be monotonic, the uniform consistency of the Free Disposal Hull (FDH) estimator introduced by Deprins {\it et al.} (1984) was shown by Korostelev {\it et al.} (1995), along with the minimax rate of uniform convergence; the uniform consistency of isotonized versions of order$-m$ frontiers introduced in Cazals {\it et al.} (2002) is proven in Daouia and Simar (2005), but rates of convergence are not available in this study. Consistency results in the $L^1$ sense were studied by Girard {\it et al.} (2005) for an estimator solving an optimization problem and by Geffroy {\it et al.} (2006) for the blockwise maxima estimator. The minimax rate of $L^1-$convergence was established by H\"ardle {\it et al.} (1995).

\vskip1ex
\noindent
Outside the field of frontier estimation, uniform convergence of the Parzen-Rosenblatt density estimator (Parzen, 1962 and Rosenblatt, 1956) was first considered by Nadaraya (1965). His results were then improved by Silverman (1978) and Stute (1982), the latter proving a law of the iterated logarithm in this context. Analogous results on kernel regression estimators were obtained by, among others, Mack and Silverman (1982), H\"ardle {\it et al.} (1988) and Einmahl and Mason (2000). The uniform consistency of isotonized versions of order$-\alpha$ quantile estimators introduced in Aragon {\it et al.} (2005) was shown in Daouia and Simar (2005). The case of estimators of left-truncated quantiles is considered in Lemdani {\it et al.} (2009). Finally, the uniform consistency of a conditional tail-index estimator is shown in Gardes and Stupfler (2013).

\vskip1ex
\noindent
The paper is organized as follows. Our main results are stated in Section~\ref{main}. The estimator is strongly uniformly consistent in a nonparametric framework. The rate of convergence is provided when the conditional survival function of $Y$ given $X=x$ belongs to the Hall class (Hall, 1982). The rate of uniform convergence is closely linked to the rate of pointwise convergence in distribution established in Girard {\it et al.} (2013). The proofs of the main results are given in Section~\ref{proof}. Auxiliary results are postponed to the Appendix. 
\section{Main results}
\label{main}
\noindent
Our results are established under the following classical condition on the kernel:

\vskip1ex

{\bf $(K)$} $K$ is a probability density function which is H\"older continuous with exponent $\eta_K$: 
$$
\exists \, c_K>0, \ \forall \, x, \, y\in \RR^d, \ |K(x)-K(y)| \leq c_K \, \| x-y \|^{\eta_K}
$$
and its support is included in $B$, the unit ball of $\RR^d$.

\vskip1ex
\noindent
Note that $(K)$ implies that $K$ is bounded with compact support. We first wish to state the uniform consistency of our estimator on a compact subset $\Omega$ of $\RR^d$ contained in the interior of $E$. To this end, three nonparametric hypotheses are introduced. The first one states the existence of the frontier $g$.

\vskip1ex

{\bf $(NP_1)$} Given $X=x$, $Y$ is positive and has a finite right endpoint $g(x)$.

\vskip1ex
\noindent 
Let $\overline{\mathcal{F}}(y \, | \, x)=\overline{F}(g(x) \, y \, | \, x)$ be the conditional survival function of the normalised random variable $Y/g(x)$ given $X=x$. The second assumption is a regularity condition on the conditional survival function of $Y$ given $X$ along the upper boundary of $S$.

\vskip1ex

{\bf $(NP_2)$} There exists $y_0\in (0, \, 1)$ such that for all $y\in [y_0, \, 1]$, $x\mapsto \overline{\mathcal{F}}(y \, | \, x)$ is continuous on $E$.

\vskip1ex
\noindent 
The third assumption, which controls the oscillation of the function $\overline{\mathcal{F}}(y\  \, | \, \cdot)$ for $y$ close to 1, can be seen as a regularity condition on the (normalised) conditional high order moment $m_{p_n}(x)/g^{p_n}(x)=\EE((Y/g(x))^{p_n} \, | \, X=x)$.

\vskip1ex

{\bf $(NP_3)$} For all $c\geq 1$,
$$
\sup_{x\in \Omega} \, \sup_{u\in B} \left| \frac{ \displaystyle\int_0^{1} y^{cp_n-1} \overline{\mathcal{F}}(y \, | \, x-h_n u) \, dy }{\displaystyle\int_0^1 y^{cp_n-1} \overline{\mathcal{F}}(y \, | \, x) \, dy } -1 \right| \rightarrow 0 \quad \textrm{ as } n\to\infty.
$$

\vskip1ex
\noindent
Let $f$ be the probability density function of $X$. The following regularity assumption is introduced:

\vskip1ex

{\bf $(A_1)$} $f$ is a positive continuous function on $E$ and $g$ is a positive H\"older continuous function on $E$ with H\"older exponent $\eta_g$.

\vskip2ex
\noindent
Before stating our first result, let us introduce some further notations. For any real-valued function $\gamma$ on $\RR^d$, the oscillation of $\gamma$ between two points $x$ and $x-h_n u$, $u\in B$, is denoted by
$$
\Delta_n^{\gamma}(x, \, u)=\gamma(x-h_n u)-\gamma(x).
$$
Finally, let $\mu_{p_n}(x)$ be the smoothed version of the conditional moment $m_{p_n}(x)$, namely 
$$
\mu_{p_n}(x)=\EE(Y^{p_n} \, K_{h_n}(x-X)) = \int_\Omega K_{h_n}(x-t) \, m_{p_n}(t) \, f(t) \, dt.
$$
Our uniform consistency result may now be stated:
\begin{Theo}
\label{uconsist}
Assume that $(NP_1-NP_3)$, $(K)$ and $(A_1)$ hold. If $p_n\to\infty$, $\dfrac{n \, h_n^d}{\log n} \displaystyle\inf_{x\in \Omega} \dfrac{\mu_{(a+1)p_n}(x)}{g^{(a+1)p_n}(x)} \to\infty$ and $p_n \, h_n^{\eta_g} \to 0$ as $n\to \infty$, then
$$
\sup_{x\in \Omega}\left| \widehat{g}_n(x) - g(x) \right| \to 0 \quad \textrm{almost surely as } n\to\infty.
$$
\end{Theo}
\noindent
As far as the conditions on $(p_n)$ and $(h_n)$ are concerned, let us highlight that, under $(A_1)$ and since $\Omega$ is compact, $f$ is uniformly continuous on $\Omega$ and $\displaystyle\inf_{\Omega} f >0$. Besides, Lemma~\ref{lemtop} implies that for $n$ large enough the ball $B(x, \, h_n)$ with center $x$ and radius $h_n$ in $\RR^d$ is contained in $E$ for every $x\in \Omega$. The uniform relative oscillation of $f$ can then be controlled as
\begin{equation}
\label{a1}
\sup_{x\in \Omega} \, \sup_{u\in B} \left| \frac{f(x-h_n u)}{f(x)}-1 \right| = \sup_{x\in \Omega} \, \sup_{u\in B} \left| \frac{\Delta_n^f(x, \, u)}{f(x)} \right| \to 0.
\end{equation}
Similarly, $\displaystyle\inf_{\Omega} g >0$ and we thus have 
\begin{equation}
\label{myeq}
\sup_{x\in \Omega} \, \sup_{u\in B} \left| \frac{\Delta_n^g(x, \, u)}{g(x)} \right| = \operatorname{O} \left( h_n^{\eta_g} \right) \to 0.
\end{equation}
Remarking that 
$$
\log\left[ \frac{g^{p_n}(x-h_n u)}{g^{p_n}(x)} \right]=p_n \log\left[ 1+\frac{\Delta_n^g(x, \, u)}{g(x)} \right]
$$
entails, if $p_n \, h_n^{\eta_g} \to 0$,
\begin{equation}
\label{a2}
\sup_{x\in \Omega} \, \sup_{u\in B} \left| \frac{g^{p_n}(x-h_n u)}{g^{p_n}(x)}-1 \right| =\operatorname{O} \left( p_n \, h_n^{\eta_g} \right).
\end{equation}
As a conclusion, the condition $p_n \, h_n^{\eta_g} \to 0$ thus makes it possible to control the oscillation of $g^{p_n}$ around $x$, uniformly in $x\in \Omega$. This condition was already introduced in Girard and Jacob (2008, 2009) and in Girard {\it et al.} (2013). 

\vskip1ex
\noindent
To give a better understanding of the conditions of Theorem~\ref{uconsist}, we introduce the semiparametric framework

\vskip1ex

{\bf $(SP)$} For all $y \in [0, \, 1]$, $\overline{\mathcal F}(y \, | \, x)=(1-y)^{\alpha(x)} \, L\left( x, \, (1-y)^{-1} \right)$, where $L$ is bounded on $\Omega \times [1, \, \infty)$ and satisfies
$$
\forall \, x\in E, \ \forall \, z\geq 1, \ L(x, \, z)=C(x)+D(x, \, z) \, z^{-\beta(x)} 
$$
where $\alpha, \, \beta$ and $C$ are positive Borel functions and $D$ is a bounded Borel function on $\Omega \times [1, \, \infty)$. 

\vskip1ex
\noindent
In model $(SP)$, the function $L(x, \, \cdot)$ is slowly varying at infinity for all $x\in E$ (see for example Bingham {\it et al.}, 1987) and belongs to the Hall class (Hall, 1982). Let us emphasize that $\alpha(x)$ drives the behavior of the distribution function of $Y$ given $X=x$ in the neighborhood of its endpoint $g(x)$. In the general context of extreme-value theory (see for instance Embrechts {\it et al.}, 1997), the conditional distribution of $Y$ given $X=x$ is said to belong to the Weibull max-domain of attraction with conditional extreme-value index $-1/\alpha(x)$. Model $(SP)$ is clearly more general than the one in Girard {\it et al.} (2013), which is restricted to the constant case $L\equiv 1$. We introduce the additional regularity condition

\vskip1ex

{\bf $(A_2)$} $\alpha$ is a H\"older continuous function on $E$ with H\"older exponent $\eta_{\alpha}$; $\beta$ and $C$ are continuous functions on $E$ and there exists $z_0\in [1, \, \infty)$ such that for all $z\geq z_0$, the map $x\mapsto D(x, \, z)$ is continuous on $E$.

\vskip1ex
\noindent
Note that if $(A_2)$ holds, 
$$
\overline{\alpha} := \displaystyle\max_{\Omega} \alpha<\infty 
$$
because $\Omega$ is compact. Our next result shows that Theorem~\ref{uconsist} holds in the semiparametric setting $(SP)$.

\begin{Coro}
\label{corSP}
Assume that $(SP)$, $(K)$ and $(A_1-A_2)$ hold. If $p_n\to\infty$, $n \, p_n^{-\overline{\alpha}} \, h_n^d/\log n \to\infty$ and $p_n \, h_n^{\eta_g}\to 0$ as $n\to \infty$, then
$$
\sup_{x\in \Omega}\left| \widehat{g}_n(x) - g(x) \right| \to 0 \quad \textrm{almost surely as } n\to\infty.
$$
\end{Coro}
\noindent
Note -- see the proof of Corollary~\ref{corSP} -- that if $(SP)$, $(K)$ and $(A_1-A_2)$ hold and $p_n \, h_n^{\eta_g}\to 0$ as $n\to \infty$, then hypothesis $(NP_3)$ holds as well. This hypothesis can therefore be considered not only as a regularity condition on the conditional high order moment $m_{p_n}(x)$ but also as a condition comparing the rates of convergence of $(1/p_n)$ and $(h_n)$ to 0.

\vskip1ex
\noindent
Our second aim is to compute the rate of convergence of the estimator~(\ref{defest}). Under hypothesis $(A_2)$, we can introduce the quantity
$$
\underline{\beta} := \displaystyle\min_{\Omega} \beta>0.
$$
Letting $w_n=\sqrt{n \, p_n^{-\overline{\alpha}+2} \, h_n^d/\log n}$, we can now state our result on the rate of uniform convergence in the semiparametric framework $(SP)$:
\begin{Theo}
\label{urateconsist}
Assume that $(SP)$, $(K)$ and $(A_1-A_2)$ hold. If $p_n\to\infty$ and
\begin{itemize}
\item ${n \, p_n^{-\overline{\alpha}} \, h_n^d}/{\log n} \to\infty$ as $n\to \infty$, 
\item $\displaystyle\limsup_{n\to\infty} w_n \, \left\{ h_n^{\eta_g} \lor p_n^{-1} \, h_n^{\eta_{\alpha}} \lor p_n^{-\underline{\beta}-1} \right\} <\infty$,
\end{itemize}
then
$$
w_n \sup_{x\in \Omega}\left| \widehat{g}_n(x) - g(x) \right| = \operatorname{O}\left( 1 \right) \quad \textrm{almost surely as } n\to\infty.
$$
\end{Theo}
\noindent
Let us highlight that the condition $n \, p_n^{-\overline{\alpha}} \, h_n^d/\log n \to\infty$ was already introduced in Corollary~\ref{corSP}. The second condition controls the bias of the estimator $\widehat{g}_n$. The term $h_n^{\eta_g}$ corresponds to the bias introduced by using a kernel smoothing, while the presence of both other terms is due to the particular structure of the semiparametric model $(SP)$. Moreover, as pointed out in Theorem~3 in Girard {\it et al.} (2013), the rate of pointwise convergence of $\widehat{g}_n(x)$ to $g(x)$ is $\sqrt{n \, p_n^{-\alpha(x)+2} \, h_n^d}$. Up to the factor $\sqrt{\log n}$, the rate of uniform convergence of $\widehat{g}_n$ is therefore the infimum (over $\Omega$) of the rate of pointwise convergence of $\widehat{g}_n(x)$ to $g(x)$.

\vskip2ex
\noindent
Theorem~\ref{urateconsist} allows us to compute the optimal rate of convergence of $\widehat{g}_n$. For the sake of simplicity, we shall consider the case when $\alpha$ is more regular than $g$ ({\it i.e.} $\eta_{\alpha}\geq \eta_g$) and $\overline{\mathcal{F}}(y\, | \, x)=(1-y)^{\alpha(x)}$ for all $y\in [0, \, 1]$ (namely, $D$ is identically zero). In that case, the conditions on $(p_n)$ and $(h_n)$ reduce to 
$$
\frac{n \, p_n^{-\overline{\alpha}} \, h_n^d}{\log n} \to\infty \quad \textrm{as } n\to \infty \quad \textrm{and} \quad \limsup_{n\to\infty} \frac{n \, p_n^{-\overline{\alpha}+2} \, h_n^{d+2\eta_g}}{\log n} <\infty.
$$
Up to the factor $\sqrt{\log n}$, the optimal rate of convergence is obtained if $p_n$ has order $n^{c_1}$ and $h_n$ has order $n^{-c_2}$, where $(c_1, \, c_2)$ is a solution of the constrained optimization problem
\begin{eqnarray*}
(c_1, \, c_2) &=& \operatorname*{argmax}_{(c, \, c') \in \Delta} \ 1+(2-\overline{\alpha})c-d c' \\
\textrm{with} \quad \Delta &=& \{ (c, \, c')\in \RR^2 \, | \, 1-\overline{\alpha} \, c - d c'\geq 0, \ 1+(2-\overline{\alpha}) c - (d+2\eta_g) c'\leq 0, \ c, \, c'> 0 \}.
\end{eqnarray*}
This yields $c_1=\eta_g/(d+\overline{\alpha} \, \eta_g)$ and $c_2=1/(d+\overline{\alpha} \, \eta_g)$, in which case the (optimal) rate of convergence has order $n^{\eta_g/(d+\overline{\alpha} \, \eta_g)}$. Let us note that this rate of convergence has been shown to be minimax by H\"ardle {\it et al.} (1995) for a particular class of densities in the case $d=1$ with a $L^1$ risk.
\section{Proofs of the main results}
\label{proof}
\noindent
Before proceeding to the proofs of our main results, we point out that, due to our hypotheses, all our results and lemmas on the behavior of $m_{p_n}(x)$, $\mu_{p_n}(x)$ and $\widehat{\mu}_{p_n}(x)$ hold as well when $p_n$ is replaced by $cp_n$, $c>1$.

\noindent 
The key idea to show Theorem~\ref{uconsist} is to prove a uniform law of large numbers for $\widehat{\mu}_{p_n}(x)$ in the nonparametric setting.
\begin{Prop}
\label{uratio}
Assume that $(NP_1-NP_3)$, $(K)$ and $(A_1)$ hold. Let $v_n=\sqrt{\dfrac{n \, h_n^d}{\log n} \displaystyle\inf_{x\in \Omega} \dfrac{\mu_{p_n}(x)}{g^{p_n}(x)}}$. If $p_n\to\infty$, $v_n\to\infty$ and $p_n \, h_n^{\eta_g} \to 0$ as $n\to \infty$, then there exists a positive constant $c>0$ such that for every $\varepsilon >0$ and every sequence of positive numbers $(\delta_n)$ converging to 0 such that $\delta_n \, v_n\to\infty$, there exists a positive constant $c_{\varepsilon}$ with
$$
\PP\left( \delta_n \, v_n \, \sup_{x\in \Omega}\left| \frac{\widehat{\mu}_{p_n}(x)}{\mu_{p_n}(x)} -1 \right| > \varepsilon \right) = \operatorname{O} \left( n^c \exp\left[ -c_{\varepsilon} \, \frac{\log n}{\delta_n^2} \right] \right).
$$
Consequently,
$$
\delta_n \, v_n \, \sup_{x\in \Omega}\left| \frac{\widehat{\mu}_{p_n}(x)}{\mu_{p_n}(x)} -1 \right| \to 0 \quad \textrm{almost surely as } n\to\infty.
$$
\end{Prop}
\noindent
{\bf Proof of Proposition~\ref{uratio}.} The proof is based on that of Lemma~1 in H\"ardle and Marron (1985). Since $\Omega$ is a compact subset of $\RR^d$, we may, for all $n\in \NN\setminus \{ 0 \}$, find a finite subset $\Omega_n$ of $\Omega$ such that: 
\begin{equation*}
\forall \, x\in \Omega, \ \exists \, \chi(x)\in \Omega_n, \ \| x-\chi(x) \| \leq n^{-\eta} \quad \textrm{and} \quad \exists \, c>0, \ \left| \Omega_n \right| = \operatorname{O}\left( n^c \right),
\end{equation*}
where $|\Omega_n|$ stands for the cardinality of $\Omega_n$, and $\eta=d^{-1}+\eta_K^{-1}$. Notice that, since $n\, h_n^d\to\infty$, one can assume that eventually $\chi(x)\in B(x, \, h_n)$ for all $x\in\Omega$. Besides, since $h_n\to 0$, we can use Lemma~\ref{lemtop} and pick $n$ so large that $B(x, \, 2h_n)\subset E$ for all $x\in\Omega$. Picking $\varepsilon>0$, and letting
\begin{eqnarray*}
T_{1, \, n} &:=& \PP\left( \delta_n \, v_n \, \sup_{x\in \Omega} \left| \frac{\widehat{\mu}_{p_n}(x)}{\mu_{p_n}(x)} -\frac{\widehat{\mu}_{p_n}(\chi(x))}{\mu_{p_n}(\chi(x))} \right| > \frac{\varepsilon}{2} \right) \\
\textrm{and} \quad T_{2, \, n} &:=& \sum_{\omega\in \Omega_n} \PP\left( \delta_n \, v_n \left| \frac{\widehat{\mu}_{p_n}(\omega)}{\mu_{p_n}(\omega)} -1 \right| > \frac{\varepsilon}{2} \right),
\end{eqnarray*}
the triangle inequality then yields
\begin{equation*}
\PP\left( \delta_n \, v_n \, \sup_{x\in \Omega}\left| \frac{\widehat{\mu}_{p_n}(x)}{\mu_{p_n}(x)} -1 \right| > \varepsilon \right) \leq T_{1, \, n} + T_{2, \, n}.
\end{equation*}
The goal of the proof is to show that 
\begin{equation*}
T_{1, \, n} + T_{2, \, n} = \operatorname{O} \left( n^c \exp\left[ -c_{\varepsilon} \, \frac{\log n}{\delta_n^2} \right] \right).
\end{equation*}
We start by controlling $T_{1, \, n}$. For all $x\in \Omega$,  
\begin{equation*}
\left| \frac{\widehat{\mu}_{p_n}(x)}{\mu_{p_n}(x)} -\frac{\widehat{\mu}_{p_n}(\chi(x))}{\mu_{p_n}(\chi(x))} \right| \leq \frac{1}{n} \sum_{i=1}^n Y_i^{p_n} \left| \frac{K_{h_n}(x-X_i)}{\mu_{p_n}(x)} - \frac{K_{h_n}(\chi(x)-X_i)}{\mu_{p_n}(\chi(x))} \right|,
\end{equation*}
and the triangle inequality entails
\begin{eqnarray*}
\left| \frac{K_{h_n}(x-X_i)}{\mu_{p_n}(x)} - \frac{K_{h_n}(\chi(x)-X_i)}{\mu_{p_n}(\chi(x))} \right| & \leq & \frac{|K_{h_n}(x-X_i)-K_{h_n}(\chi(x)-X_i)|}{\mu_{p_n}(x)} \\
  & + &  \frac{|\mu_{p_n}(x)- \mu_{p_n}(\chi(x))|}{\mu_{p_n}(x) \, \mu_{p_n}(\chi(x))} \, K_{h_n}(\chi(x)-X_i).
\end{eqnarray*}
Using hypothesis $(K)$ and Lemma~\ref{lemosc}, there exists a positive constant $\kappa$ such that, for $n$ large enough,
\begin{equation*}
\sup_{x\in \Omega} \left\{ \mu_{p_n}(x) \left| \frac{K_{h_n}(x-X_i)}{\mu_{p_n}(x)} - \frac{K_{h_n}(\chi(x)-X_i)}{\mu_{p_n}(\chi(x))} \right| \right\} \leq \frac{\kappa}{h_n^d} \left[ \frac{n^{-\eta}}{h_n} \right]^{\eta_K} \ind_{\{ X\in B(x, \, h_n) \cup B(\chi(x), \, h_n) \}}. 
\end{equation*}
Since the support of the random variable $K_{h_n}(\chi(x)-X_i)$ is included in $B(x, \, 2h_n)$, one has
\begin{equation*}
\sup_{x\in \Omega} \left| \frac{\widehat{\mu}_{p_n}(x)}{\mu_{p_n}(x)} -\frac{\widehat{\mu}_{p_n}(\chi(x))}{\mu_{p_n}(\chi(x))} \right| \leq \kappa \left[ \frac{n^{-\eta}}{h_n} \right]^{\eta_K} \sup_{x\in \Omega} \frac{1}{\mu_{p_n}(x)} \left| \frac{1}{n \, h_n^d} \sum_{i=1}^n Y_i^{p_n} \ind_{\{ X_i\in B(x, \, 2h_n) \}} \right|.
\end{equation*}
For all $x\in\Omega$,
\begin{equation*}
\frac{1}{n} \sum_{i=1}^n Y_i^{p_n} \ind_{\{ X_i\in B(x, \, 2h_n) \}} \leq \sup_{B(x, \, 2h_n)} g^{p_n}
\end{equation*}
almost surely, and in view of~(\ref{a2}), it follows that
\begin{equation*}
\sup_{x\in \Omega} \left| \frac{\widehat{\mu}_{p_n}(x)}{\mu_{p_n}(x)} -\frac{\widehat{\mu}_{p_n}(\chi(x))}{\mu_{p_n}(\chi(x))} \right| \leq 2\kappa \left[ \frac{n^{-\eta}}{h_n} \right]^{\eta_K} \frac{1}{h_n^d} \sup_{x\in \Omega} \frac{g^{p_n}(x)}{\mu_{p_n}(x)}
\end{equation*}
for $n$ large enough. Finally, $n \, h_n^d\to\infty$ implies
\begin{equation*}
\left[ \frac{n^{-\eta}}{h_n} \right]^{\eta_K} = \left[ \frac{1}{n \, h_n^d} \right]^{\eta_K/d} \, \frac{1}{n} = \operatorname{o}\left( \frac{1}{n} \right)
\end{equation*}
and therefore, we have the following bound: 
\begin{equation*}
\delta_n \, v_n \, \sup_{x\in \Omega} \left| \frac{\widehat{\mu}_{p_n}(x)}{\mu_{p_n}(x)} -\frac{\widehat{\mu}_{p_n}(\chi(x))}{\mu_{p_n}(\chi(x))} \right| \leq 2\kappa \frac{\delta_n}{v_n \, \log n} \to 0
\end{equation*}
as $n\to\infty$. Hence $T_{1, \, n}=0$ eventually.

\vskip1ex
\noindent
Let us now control $T_{2, \, n}$. To this end, pick $\omega\in \Omega_n$ and introduce
\begin{equation*}
Z_{n, \, i}(\omega)=\dfrac{Y_i^{p_n}}{\displaystyle\sup_{B(\omega, \, h_n)} g^{p_n}} \, K \left( \dfrac{\omega-X_i}{h_n} \right).
\end{equation*}
Remark that $\left| Z_{n, \, i}(\omega)-\EE(Z_{n, \, i}(\omega)) \right| \leq \displaystyle\sup_{B} K$ almost surely and thus
\begin{equation*}
h_n^d \, \frac{\widehat{\mu}_{p_n}(\omega)-\mu_{p_n}(\omega)}{\displaystyle\sup_{B(\omega, \, h_n)} g^{p_n}} = \frac{1}{n} \sum_{i=1}^n \big\{ Z_{n, \, i}(\omega)-\EE(Z_{n, \, i}(\omega)) \big\}
\end{equation*}
is a mean of bounded, centered, independent and identically distributed random variables. Defining 
\begin{eqnarray*}
\tau_n(\omega) &:=& \frac{\varepsilon}{2 \displaystyle\sup_B K} \, \frac{1}{\delta_n \, v_n} \, \frac{n\, \mu_{p_n}(\omega) \, h_n^d}{\displaystyle\sup_{B(\omega, \, h_n)} g^{p_n}} \\ 
\textrm{and} \quad \lambda_n(\omega) &:=& \frac{\varepsilon}{2} \, \sup_B K \, \frac{1}{\delta_n \, v_n} \, \frac{\mu_{p_n}(\omega) \, h_n^d}{\displaystyle\sup_{B(\omega, \, h_n)} g^{p_n}} \, \frac{1}{\Var(Z_{n, \, 1}(\omega))},
\end{eqnarray*}
Bernstein's inequality (see Hoeffding, 1963) yields, for all $\varepsilon>0$, 
\begin{eqnarray*}
\PP\left( \delta_n \, v_n \left| \frac{\widehat{\mu}_{p_n}(\omega)}{\mu_{p_n}(\omega)} -1 \right| > \frac{\varepsilon}{2} \right) &=& \PP\left( h_n^d \left| \frac{\widehat{\mu}_{p_n}(\omega)-\mu_{p_n}(\omega)}{\displaystyle\sup_{B(\omega, \, h_n)} g^{p_n}} \right| > \frac{\varepsilon}{2} \, \frac{1}{\delta_n \, v_n} \, \frac{\mu_{p_n}(\omega) \, h_n^d}{\displaystyle\sup_{B(\omega, \, h_n)} g^{p_n}} \right) \\ 
 & \leq & 2\exp\left( -\frac{\tau_n(\omega) \lambda_n(\omega)}{2(1+\lambda_n(\omega)/3)} \right).
\end{eqnarray*}
Using once again~(\ref{a2}), we get, for $n$ large enough,
\begin{equation*}
\inf_{\omega\in \Omega_n} \tau_n(\omega) \geq \frac{\varepsilon}{4 \displaystyle\sup_B K} \, \frac{v_n \, \log n}{\delta_n}.
\end{equation*}
Moreover, for all $\omega\in\Omega_n$, 
\begin{equation*}
\frac{1}{\lambda_n(\omega)} = \frac{2}{\varepsilon \displaystyle\sup_B K} \, \delta_n \, v_n \, \sup_{B(\omega, \, h_n)} g^{p_n} \, h_n^{-d} \left[ \frac{\EE(Z_{n, \, 1}^2(\omega))}{\mu_{p_n}(\omega)} - \frac{\left[ \EE(Z_{n, \, 1}(\omega)) \right]^2}{\mu_{p_n}(\omega)} \right],
\end{equation*}
and since $\displaystyle\sup_{B(\omega, \, h_n)} g^{p_n} \, h_n^{-d} \, Z_{n, \, 1}(\omega) = Y_1^{p_n} \, K_{h_n}(\omega-X_1)$, it follows that
\begin{equation*}
\sup_{B(\omega, \, h_n)} g^{p_n} \, h_n^{-d} \left[ \frac{\EE(Z_{n, \, 1}^2(\omega))}{\mu_{p_n}(\omega)} - \frac{\left[ \EE(Z_{n, \, 1}(\omega)) \right]^2}{\mu_{p_n}(\omega)} \right] \leq \sup_B K ,
\end{equation*}
so that 
\begin{equation*}
\sup_{\omega\in \Omega_n} \frac{1}{\lambda_n(\omega)} \leq \frac{2}{\varepsilon} \, \delta_n \, v_n.
\end{equation*}
Remarking that the function $x\mapsto 1/[2(x+1/3)]$ is decreasing on $\RR_+$, there exists a constant $c_{\varepsilon}>0$ such that, for all $\omega\in \Omega_n$,
\begin{equation*}
 \PP\left( \delta_n \, v_n \left| \frac{\widehat{\mu}_{p_n}(\omega)}{\mu_{p_n}(\omega)} -1 \right| > \frac{\varepsilon}{2} \right) \leq 2\exp\left( -c_{\varepsilon} \, \frac{\log n}{\delta_n^2} \right),
\end{equation*}
for all $n$ large enough. Taking into account that $|\Omega_n|=\operatorname{O}(n^c)$, this implies that 
\begin{equation*}
T_{2, \, n} = \operatorname{O} \left( n^c \exp\left[ -c_{\varepsilon} \, \frac{\log n}{\delta_n^2} \right] \right).
\end{equation*}
Notice now that the above bound yields
\begin{equation*}
\forall \, \varepsilon>0, \ \sum_n \PP\left( \delta_n \, v_n \, \sup_{x\in \Omega}\left| \frac{\widehat{\mu}_{p_n}(x)}{\mu_{p_n}(x)} -1 \right| > \varepsilon \right) <\infty
\end{equation*}
and use Borel-Cantelli's lemma to get the final part of the result. \CQFD

\vskip2ex
\noindent
With Proposition~\ref{uratio} at hand, we can now prove Theorem~\ref{uconsist}.

\vskip2ex
\noindent
{\bf Proof of Theorem~\ref{uconsist}.} Since $g$ is positive and continuous on the compact set $\Omega$, it is bounded from below by a positive constant. It is then enough to prove that 
\begin{equation*}
\sup_{x\in \Omega}\left| \frac{1}{\widehat{g}_n(x)}- \frac{1}{g(x)} \right| \to 0 \quad \textrm{almost surely as } n\to\infty.
\end{equation*}
To this end, notice that 
\begin{eqnarray*}
\frac{\widehat{\mu}_{(a+1)p_n}(x)}{\widehat{\mu}_{(a+1)p_n+1}(x)} &=& \frac{\mu_{(a+1)p_n}(x)}{\mu_{(a+1)p_n+1}(x)} \, \frac{\widehat{\mu}_{(a+1)p_n}(x)}{\mu_{(a+1)p_n}(x)} \, \left[ \frac{\widehat{\mu}_{(a+1)p_n+1}(x)}{\mu_{(a+1)p_n+1}(x)} \right]^{-1} \\
 \textrm{and} \quad \frac{\widehat{\mu}_{p_n}(x)}{\widehat{\mu}_{p_n+1}(x)} &=& \frac{\mu_{p_n}(x)}{\mu_{p_n+1}(x)} \, \frac{\widehat{\mu}_{p_n}(x)}{\mu_{p_n}(x)} \, \left[ \frac{\widehat{\mu}_{p_n+1}(x)}{\mu_{p_n+1}(x)} \right]^{-1}.
\end{eqnarray*}
Using again the positivity and the continuity of $g$ on the compact set $\Omega$, Lemma~\ref{lemuratiodet}(iii) yields 
\begin{equation*}
\sup_{x\in \Omega}\left| \frac{\mu_{p_n+1}(x)}{\mu_{p_n}(x)} -g(x) \right| \to 0 \quad \textrm{and} \quad \sup_{x\in \Omega}\left| \frac{\mu_{(a+1)p_n+1}(x)}{\mu_{(a+1)p_n}(x)} -g(x) \right| \to 0.
\end{equation*}
Since $\mu_{(a+1)p_n}(x)/g^{(a+1)p_n}(x) \leq \mu_{p_n}(x)/g^{p_n}(x) \, (1+\operatorname{o}(1))$ uniformly in $x\in\Omega$, Proposition~\ref{uratio} entails
\begin{equation}
\label{e1}
\sup_{x\in \Omega} \left| \frac{\widehat{\mu}_{(a+1)p_n}(x)}{\widehat{\mu}_{(a+1)p_n+1}(x)} - \frac{1}{g(x)} \right| \to 0 \quad \textrm{and} \quad \sup_{x\in \Omega} \left| \frac{\widehat{\mu}_{p_n}(x)}{\widehat{\mu}_{p_n+1}(x)} - \frac{1}{g(x)} \right| \to 0
\end{equation}
almost surely as $n\to\infty$. The result follows by reporting~(\ref{e1}) into~(\ref{defest}). \CQFD

\vskip2ex
\noindent
Before proving Corollary~\ref{corSP}, a further examination of the behavior of the high order moment $\mu_{p_n}(x)$ is needed. The next result gives a uniform equivalent of the moment $\mu_{p_n}(x)$ in the semiparametric framework.
\begin{Prop}
\label{propequivmu}
Assume that $(SP)$, $(K)$, $(A_1-A_2)$ hold, $p_n\to\infty$ and $p_n \, h_n^{\eta_g} \to 0$ as $n\to \infty$. Then
\begin{equation*}
\sup_{x\in\Omega} \left| \frac{\mu_{p_n}(x)}{f(x) \, C(x) \, \Gamma(\alpha(x)+1) \, g^{p_n}(x) \, p_n^{-\alpha(x)}} - 1 \right| \to 0 \quad \textrm{as } n\to\infty.
\end{equation*}
\end{Prop}
\noindent
{\bf Proof of Proposition~\ref{propequivmu}.} Let us introduce $\overline{\mathcal{F}}_{\gamma}(y \, | \, x)=(1-y)^{\gamma(x)}$ for all $y\in [0, \, 1]$. In the semiparametric setting $(SP)$, $\overline{\mathcal{F}}(\cdot \, | \, x)$ can be written as
\begin{equation*}
\forall \, y\in [0, \, 1], \ \overline{\mathcal{F}}(y \, | \, x)=C(x) \, \overline{\mathcal{F}}_{\alpha}(y \, | \, x) + D\left( x, \, (1-y)^{-1} \right) \, \overline{\mathcal{F}}_{\alpha+\beta}(y \, | \, x). 
\end{equation*}
Using Lemma~\ref{lemtop}, we can pick $n$ large enough such that $B(x, \, h_n)\subset E$ for all $x\in\Omega$. Pick then $x\in \Omega$, and set
\begin{eqnarray}
\label{notation_M}
M_n(p_n, \, x) &:=& \int_E f(v) \, C(v) \, g^{p_n}(v) \, K_{h_n}(x-v) \left[ p_n \int_{0}^{\infty} y^{p_n-1} \overline{\mathcal{F}}_{\alpha}(y \, | \, v) \, dy \right] dv \\
\nonumber
			   &=& \int_B (f C g^{p_n})(x-h_n u) \, p_n \, b(p_n, \, \alpha(x-h_n u)+1) \, K(u) \, du
\end{eqnarray}
where $b(x, \, y)=\displaystyle\int_0^1 t^{x-1} \, (1-t)^{y-1} \, dt$ is the Beta function. With these notations, the high order moment $\mu_{p_n}(x)$ can be rewritten as
\begin{equation}
\label{devmu}
\mu_{p_n}(x)=M_n(p_n, \, x)[1+\varepsilon_n(p_n, \, x)] \quad \textrm{where} \quad \varepsilon_n(p_n, \, x)=\frac{E_n(p_n, \, x)}{M_n(p_n, \, x)}
\end{equation}
and with
\begin{eqnarray}
\label{notation_E}
E_n(p_n, \, x) &:=& \int_B (f g^{p_n})(x-h_n u) \, p_n \, \mathcal{I}_{\alpha+\beta, \, D}(p_n, \, x-h_n u) \, K(u) \, du \\
\label{notation_I}
\textrm{where } \ \mathcal{I}_{\alpha+\beta, \, D}(p_n, \, v) &:=& \int_0^1 y^{p_n-1} \, \overline{\mathcal{F}}_{\alpha+\beta}(y \, | \, v) \, D\left(v, \, (1-y)^{-1} \right) \, dy.
\end{eqnarray}
Lemma~\ref{lemerr} and~(\ref{devmu}) entail 
\begin{equation*}
\sup_{x\in\Omega} \left| \frac{\mu_{p_n}(x)}{M_n(p_n, \, x)} - 1 \right| \to 0 \quad \textrm{as } n\to\infty.
\end{equation*}
It is therefore enough to show that 
\begin{equation*}
\sup_{x\in\Omega} \left| \frac{M_n(p_n, \, x)}{f(x) \, C(x) \, \Gamma(\alpha(x)+1) \, g^{p_n}(x) \, p_n^{-\alpha(x)}} - 1 \right| \to 0 \quad \textrm{as } n\to\infty.
\end{equation*}
Lemma~\ref{lemdvpt} establishes that
\begin{equation*}
\sup_{x\in\Omega} \left| \dfrac{M_n(p_n, \, x)}{ f(x) \, C(x) \, \alpha(x) \, g^{p_n}(x) \, b(p_n+1, \, \alpha(x))} - 1 \right| \to 0 \quad \textrm{as } n\to\infty.
\end{equation*}
Finally, Lemma~\ref{lemstirunif} gives
\begin{equation*}
\sup_{x\in \Omega} \left| \frac{\alpha(x) \, b(p_n+1, \, \alpha(x))}{\Gamma(\alpha(x)+1) \, p_n^{-\alpha(x)}} -1 \right| \to 0 \quad \textrm{as } n\to\infty
\end{equation*}
and the result is proven. \CQFD

\vskip2ex
\noindent
Corollary~\ref{corSP} can now be shown. 

\vskip2ex
\noindent
{\bf Proof of Corollary~\ref{corSP}.} It is enough to check that the hypotheses of Theorem~\ref{uconsist} are satisfied. This is clearly the case for $(NP_1)$ and $(NP_2)$; besides, for all $c\geq 1$, Proposition~\ref{propequivmu} yields
$$
\sup_{x\in\Omega} \left| \frac{\mu_{cp_n}(x)}{f(x) \, C(x) \, \Gamma(\alpha(x)+1) \, g^{cp_n}(x) \, (cp_n)^{-\alpha(x)}} - 1 \right| \to 0 \quad \textrm{as } n\to\infty.
$$
Using Lemma~\ref{lemuratiodet} then gives
$$
\sup_{x\in\Omega} \left| \frac{\displaystyle\int_0^1 y^{cp_n-1} \overline{\mathcal{F}}(y \, | \, x) \, dy }{C(x) \, \Gamma(\alpha(x)+1) \, \, (cp_n)^{-\alpha(x)}} - 1 \right| \to 0 \quad \textrm{as } n\to\infty.
$$
The hypothesis $p_n \, h_n^{\eta_g} \to 0$ thus makes it clear that $(NP_3)$ holds as well in this setting. Finally, Proposition~\ref{propequivmu} entails 
$$
\sup_{x\in\Omega} \left| \frac{\mu_{(a+1)p_n}(x)/g^{(a+1)p_n}(x)}{f(x) \, C(x) \, \Gamma(\alpha(x)+1) \, [(a+1)p_n]^{-\alpha(x)}} - 1 \right| \to 0 \quad \textrm{as } n\to\infty.
$$
Consequently, for $n$ large enough there exists some positive constant $\varepsilon>0$ such that 
$$
\frac{n \, h_n^d}{\log n} \inf_{x\in \Omega} \frac{\mu_{(a+1)p_n}(x)}{g^{(a+1)p_n}(x)} \geq \varepsilon \frac{n \, p_n^{-\overline{\alpha}} \, h_n^d}{\log n} \to\infty \ \mbox{ as } \ n\to\infty
$$
which concludes the proof. \CQFD

\vskip2ex
\noindent
In order to prove Theorem~\ref{urateconsist}, since the expression of our frontier estimator involve ratios such as $\widehat{\mu}_{p_n}(x)/\widehat{\mu}_{p_n+1}(x)$, we shall first compute an asymptotic expansion of $\mu_{p_n}(x)/\mu_{p_n+1}(x)$:
\begin{Prop}
\label{propdeux}
Assume that $(SP)$, $(K)$ and $(A_1-A_2)$ hold. If $p_n\to\infty$ and $p_n \, h_n^{\eta_g} \to 0$, then 
\begin{equation*}
\sup_{x\in\Omega} \left\{ \frac{1}{h_n^{\eta_g} \lor p_n^{-1} \, h_n^{\eta_{\alpha}} \lor p_n^{-\beta(x)-1}} \left| \frac{\mu_{p_n}(x)}{\mu_{p_n+1}(x)} - \frac{1}{g(x)} \, \left[ 1+\frac{\alpha(x)}{p_n+1} \right] \right| \right\} = \operatorname{O}(1).
\end{equation*}
\end{Prop}
\noindent
{\bf Proof of Proposition~\ref{propdeux}.} Remark that, with the notations of Proposition~\ref{propequivmu} above, we have
\begin{eqnarray}
\label{errmu}
\frac{\mu_{p_n}(x)}{\mu_{p_n+1}(x)} &=& \frac{M_n(p_n, \, x)}{M_n(p_n+1, \, x)} \, [1+\tau_n(p_n, \, x)] \\
\nonumber
\textrm{where} \quad \tau_n(p_n, \, x) &:=& \dfrac{\varepsilon_n(p_n, \, x)-\varepsilon_n(p_n+1, \, x)}{1+\varepsilon_n(p_n+1, \, x)}.
\end{eqnarray}
Using Lemma~\ref{lemtop}, we can pick $n$ large enough such that $B(x, \, h_n)\subset E$ for all $x\in\Omega$. Recall then the notations of Lemma~\ref{lemdvpt} and write
\begin{equation*}
\sup_{x\in\Omega} \left| \left\{ \dfrac{1}{g(x)} \left[ 1+\dfrac{\alpha(x)}{p_n+1} \right] \right\}^{-1} \dfrac{M_n(p_n, \, x)}{M_n(p_n+1, \, x)} - \frac{\displaystyle\int_B \mathcal{L}_n(p_n, \, x, \, u) \, K(u) \, du}{\displaystyle\int_B \mathcal{L}_n(p_n+1, \, x, \, u) \, K(u) \, du} \right| = \operatorname{O}\left( h_n^{\eta_g} \lor \frac{h_n^{\eta_{\alpha}}}{p_n} \right).
\end{equation*}
Since $\mathcal{L}_n(p_n+1, \, x, \, u) > 0$, it follows that
\begin{equation*}
\sup_{x\in \Omega} \left| \frac{\displaystyle\int_B \mathcal{L}_n(p_n, \, x, \, u) \, K(u) \, du}{\displaystyle\int_B \mathcal{L}_n(p_n+1, \, x, \, u) \, K(u) \, du} - 1 \right| \leq \sup_{x\in \Omega} \, \sup_{u\in B} \left| \frac{\mathcal{L}_n(p_n, \, x, \, u)}{\mathcal{L}_n(p_n+1, \, x, \, u)} - 1 \right|
= \operatorname{O}\left( h_n^{\eta_g} \lor \frac{h_n^{\eta_{\alpha}}}{p_n} \right).
\end{equation*}
Lemma~\ref{lemdvpt} entails
\begin{equation*}
\sup_{x\in\Omega} \left| \left\{ \dfrac{1}{g(x)} \left[ 1+\dfrac{\alpha(x)}{p_n+1} \right] \right\}^{-1} \dfrac{M_n(p_n, \, x)}{M_n(p_n+1, \, x)} - 1 \right| = \operatorname{O}\left( h_n^{\eta_g} \lor \frac{h_n^{\eta_{\alpha}}}{p_n} \right).
\end{equation*}
Besides, applying Lemma~\ref{lemerr} yields $\displaystyle\sup_{x\in\Omega} \left| p_n^{\beta(x)+1} \, \tau_n(p_n, \, x) \right|=\operatorname{O}(1)$. Replacing in~(\ref{errmu}) concludes the proof of Proposition~\ref{propdeux}. \CQFD

\vskip2ex
\noindent
We can now give a proof of Theorem~\ref{urateconsist}.

\vskip2ex
\noindent
{\bf Proof of Theorem~\ref{urateconsist}.} Since, by Theorem~\ref{uconsist}, $\displaystyle\sup_{x\in \Omega} |\widehat{g}_n(x) - g(x)| \to 0$ almost surely, it is enough to prove that 
\begin{equation*}
w_n\sup_{x\in \Omega}\left| \frac{1}{\widehat{g}_n(x)} - \frac{1}{g(x)} \right| = \operatorname{O}\left(1\right) \quad \textrm{almost surely as } n\to\infty.
\end{equation*}
Introducing
\begin{equation*}
\frac{1}{G_n(x)} = \frac{1}{a p_n} \left[ ((a+1)p_n+1) \frac{\mu_{(a+1)p_n}(x)}{\mu_{(a+1)p_n+1}(x)} - (p_n+1) \frac{\mu_{p_n}(x)}{\mu_{p_n+1}(x)} \right] \textrm{ and } \ \xi_n(x)=\frac{1}{\widehat{g}_n(x)} - \frac{1}{G_n(x)}
\end{equation*} 
the quantity of interest can be expanded as
\begin{equation*}
\frac{1}{\widehat{g}_n(x)} - \frac{1}{g(x)} = \xi_n(x)+ \left[ \frac{1}{G_n(x)} - \frac{1}{g(x)} \right]  .
\end{equation*}
Both terms are considered separately. The bias term is readily controlled by Proposition~\ref{propdeux}:
\begin{equation*}
w_n \sup_{x\in\Omega} \left| \frac{1}{G_n(x)} - \frac{1}{g(x)} \right| = \operatorname{O}\left( w_n \left\{ h_n^{\eta_g} \lor \frac{h_n^{\eta_{\alpha}}}{p_n} \lor p_n^{-\underline{\beta}-1} \right\} \right) = \operatorname{O}(1)
\end{equation*}
in view of the hypotheses on $(p_n)$ and $(h_n)$. Let us now consider the random term $\xi_n(x)$. Lemma~\ref{lemlin} shows that
\begin{equation*}
\xi_n(x) = \frac{1}{a p_n} \left[ \zeta_n^{(1)}(x) - \zeta_n^{(2)}(x) + \left(\frac{\mu_{p_n+1}(x)}{\widehat{\mu}_{p_n+1}(x)}-1 \right) \zeta_n^{(1)}(x) - \left( \frac{\mu_{(a+1) p_n+1}(x)}{\widehat{\mu}_{(a+1) p_n+1}(x)}-1 \right) \zeta_n^{(2)}(x) \right].
\end{equation*}
In view of Proposition~\ref{uratio}, it is therefore sufficient to show that 
\begin{equation}
\label{random}
\frac{w_n}{p_n} \, \sup_{x\in\Omega} \left| \zeta_n^{(1)}(x) \right| = \operatorname{O}\left( 1\right) \quad \textrm{and} \quad \frac{w_n}{p_n} \, \sup_{x\in\Omega} \left| \zeta_n^{(2)}(x) \right| = \operatorname{O}\left( 1 \right)
\end{equation}
almost surely as $n\to\infty$. We shall only prove the result for $\zeta_n^{(1)}(x)$, since the result will then be obtained for $\zeta_n^{(2)}(x)$ by replacing $p_n$ with $(a+1)p_n$. To this end, we mimick the proof of Proposition~\ref{uratio}. For all $n\in \NN\setminus\{ 0 \}$, let $\Omega_n$ be a finite subset of $\Omega$ such that: 
\begin{equation*}
\forall \, x\in \Omega, \ \exists \, \chi(x)\in \Omega_n, \ \| x-\chi(x) \| \leq n^{-\eta} \quad \textrm{and} \quad \exists \, c>0, \ \left| \Omega_n \right| = \operatorname{O}\left( n^c \right),
\end{equation*}
where $\eta=d^{-1}+\eta_K^{-1} \left[ 1+\overline{\alpha}^{-1} \right]$ and assume that $n$ is large enough so that $\chi(x)\in B(x, \, h_n)$ and, by Lemma~\ref{lemtop}, such that $B(x, \, 2h_n)\subset E$ for all $x\in\Omega$. Pick $\varepsilon>0$ and an arbitrary positive sequence $(\delta_n)$ converging to 0, and let 
\begin{eqnarray*}
T_{1, \, n} &:=& \PP\left( \delta_n \, \frac{w_n}{p_n} \, \sup_{x\in \Omega}\left| \zeta_n^{(1)}(x) -\zeta_n^{(1)}(\chi(x)) \right|  > \frac{\varepsilon}{2} \right) \\
\textrm{and} \quad T_{2, \, n} &:=& \sum_{\omega\in \Omega_n} \PP\left( \delta_n \, \frac{w_n}{p_n} \left| \zeta_n^{(1)}(\omega) \right|  > \frac{\varepsilon}{2} \right).
\end{eqnarray*}
The goal is then to show that both series $\sum_n T_{1, \, n}$ and $\sum_n T_{2, \, n}$ converge. Noting that 
$$
\delta_n \leq \delta_n \lor \sqrt{\frac{p_n}{w_n}}
$$
we shall assume without loss of generality that $\delta_n \sqrt{n \, p_n^{-\overline{\alpha}} \, h_n^d/\log n} \to\infty$. Let first 
\begin{eqnarray*}
T_{3, \, n} &:=& \PP\left( \delta_n \, w_n \, \sup_{x\in \Omega}\left| \frac{\mu_{p_n}(\chi(x))}{\mu_{p_n+1}(\chi(x))} \left[ \frac{\widehat{\mu}_{p_n}(x)}{\mu_{p_n}(x)} - \frac{\widehat{\mu}_{p_n}(\chi(x))}{\mu_{p_n}(\chi(x))} \right] \right| > \frac{\varepsilon}{16} \right), \\
T_{4, \, n} &:=& \PP\left( \delta_n \, w_n \, \sup_{x\in \Omega} \left| \frac{\mu_{p_n}(\chi(x))}{\mu_{p_n+1}(\chi(x))} \left[ \frac{\widehat{\mu}_{p_n+1}(x)}{\mu_{p_n+1}(x)} - \frac{\widehat{\mu}_{p_n+1}(\chi(x))}{\mu_{p_n+1}(\chi(x))} \right] \right| > \frac{\varepsilon}{16} \right), \\
T_{5, \, n} &:=& \PP\left( \delta_n \, w_n \, \sup_{x\in \Omega}\left| \left[ \frac{\mu_{p_n}(x)}{\mu_{p_n+1}(x)} - \frac{\mu_{p_n}(\chi(x))}{\mu_{p_n+1}(\chi(x))} \right] \left[ \frac{\widehat{\mu}_{p_n}(x)}{\mu_{p_n}(x)} - 1 \right] \right| > \frac{\varepsilon}{16} \right), \\
\textrm{and} \quad T_{6, \, n} &:=& \PP\left( \delta_n \, w_n \, \sup_{x\in \Omega}\left| \left[ \frac{\mu_{p_n}(x)}{\mu_{p_n+1}(x)} - \frac{\mu_{p_n}(\chi(x))}{\mu_{p_n+1}(\chi(x))} \right] \left[ \frac{\widehat{\mu}_{p_n+1}(x)}{\mu_{p_n+1}(x)} - 1 \right] \right| > \frac{\varepsilon}{16} \right), 
\end{eqnarray*}
so that for all sufficiently large $n$, $T_{1, \, n} \leq T_{3, \, n} + T_{4, \, n} + T_{5, \, n} + T_{6, \, n}$. A proof similar to the one of Proposition~\ref{uratio} gives the bound
\begin{equation*}
\sup_{x\in \Omega} \left| \frac{\widehat{\mu}_{p_n}(x)}{\mu_{p_n}(x)} -\frac{\widehat{\mu}_{p_n}(\chi(x))}{\mu_{p_n}(\chi(x))} \right| \leq \kappa \left[ \frac{n^{-\eta}}{h_n} \right]^{\eta_K} \frac{1}{h_n^d} \sup_{x\in \Omega} p_n^{\alpha(x)} 
\end{equation*}
for $n$ large enough, where $\kappa$ is a positive constant. Remark that $n \, p_n^{-\overline{\alpha}} \to \infty$ and $n\, h_n^d \to\infty$ yield
\begin{equation*}
p_n \left[ \frac{n^{-\eta}}{h_n} \right]^{\eta_K} = \left[ \frac{1}{n\, p_n^{-\overline{\alpha}}} \right]^{1/\overline{\alpha}} \, \left[ \frac{1}{n \, h_n^d} \right]^{\eta_K/d} \, \frac{1}{n} = \operatorname{o}\left( \frac{1}{n} \right).
\end{equation*}
Recalling that, from Proposition~\ref{uratio},
\begin{equation*}
v_n=\sqrt{\frac{n\, h_n^d}{\log n}\displaystyle\inf_{x\in \Omega} \dfrac{\mu_{p_n}(x)}{g^{p_n}(x)}}
\end{equation*}
Proposition~\ref{propequivmu} yields $w_n = p_n\, v_n$ and therefore, applying Proposition~\ref{propdeux}, $T_{3, \, n}=0$ and $T_{4, \, n}=0$ eventually as $n\to\infty$, so that $\sum_n T_{3, \, n}$ and $\sum_n T_{4, \, n}$ converge. Furthermore, since $\chi(x)\in B(x, \, h_n)$, Proposition~\ref{propdeux} entails
\begin{equation}
\label{diffmu}
\sup_{x\in \Omega}\left| \frac{\mu_{p_n}(x)}{\mu_{p_n+1}(x)} - \frac{\mu_{p_n}(\chi(x))}{\mu_{p_n+1}(\chi(x))} \right| = \operatorname{O}\left( h_n^{\eta_g} \lor \frac{h_n^{\eta_{\alpha}}}{p_n} \lor p_n^{-\underline{\beta}-1} \right).
\end{equation}
Using once again the equality $w_n=p_n \, v_n$ and~(\ref{diffmu}) together with Proposition~\ref{uratio} shows that $\sum_n T_{5, \, n}$ and $\sum_n T_{6, \, n}$ converge. As a consequence, $\sum_n T_{1, \, n}$ converges.

\vskip1ex
\noindent
To control $T_{2, \, n}$, we shall, as in the proof of Proposition~\ref{uratio}, show that there exists a positive constant $c_{\varepsilon}$ such that for all sufficiently large $n$,
\begin{equation*}
\forall \, \omega\in\Omega_n, \ \PP\left( \delta_n \, \frac{w_n}{p_n} \left| \zeta_n^{(1)}(\omega) \right| > \frac{\varepsilon}{2} \right) \leq \exp\left( -c_{\varepsilon} \, \frac{\log n}{\delta_n^2} \right).
\end{equation*}
Pick $\omega\in\Omega_n$ and let us consider the random variables
\begin{equation*}
S_{n, \, i}(\omega) = Y_i^{p_n} \left[ -1 + \frac{\mu_{p_n}(\omega)}{\mu_{p_n+1}(\omega)} \, Y_i \right] K_{h_n}(\omega-X_i), \ i=1, \, \ldots, \, n
\end{equation*}
such that
\begin{equation}
\label{defzeta1}
\zeta_n^{(1)}(\omega) = \frac{p_n+1}{\mu_{p_n+1}(\omega)} \, \frac{1}{n} \sum_{i=1}^n S_{n, \, i}(\omega).
\end{equation}
Let now $U_{n, \, i}(\omega)=Y_i \left/\displaystyle\sup_{B(\omega, \, h_n)} g \right.$, so that $U_{n, \, i}(\omega) \leq 1$ given $\{ X_i\in B(\omega, \, h_n) \}$. It follows that
\begin{equation*}
\frac{h_n^d}{\displaystyle\sup_{B(\omega, \, h_n)} g^{p_n}} \, S_{n, \, i}(\omega) = U_{n, \, i}^{p_n}(\omega) \left[ -1 + \sup_{B(\omega, \, h_n)} g \, \frac{\mu_{p_n}(\omega)}{\mu_{p_n+1}(\omega)} \, U_{n, \, i}(\omega) \right] K\left( \frac{\omega-X_i}{h_n} \right).
\end{equation*}
Using Proposition~\ref{propdeux}, the H\"older continuity of $g$ and the fact that $p_n \, h_n^{\eta_g} \to 0$ therefore yields, for $n$ sufficiently large, 
\begin{equation*}
(p_n+1) \sup_{\omega\in \Omega_n} \left| \frac{h_n^d}{\displaystyle\sup_{B(\omega, \, h_n)} g^{p_n}} \, S_{n, \, i}(\omega) - U_{n, \, i}^{p_n}(\omega) \left[ U_{n, \, i}(\omega)-1 \right] K\left( \frac{\omega-X_i}{h_n} \right) \right| \leq \kappa'
\end{equation*}
where $\kappa'$ is a positive constant. Some straightforward real analysis shows that
\begin{equation*}
(p_n+1) \sup_{u\in[0, \, 1]} u^{p_n} (1-u) = \left[ 1-\frac{1}{p_n+1} \right]^{p_n} \to e^{-1} <\infty.
\end{equation*}
Consequently, there exists a positive constant $\kappa''$ such that, for $n$ large enough,
\begin{equation*}
(p_n+1) \sup_{\omega\in \Omega_n} \left| \frac{h_n^d}{\displaystyle\sup_{B(\omega, \, h_n)} g^{p_n}} S_{n, \, i}(\omega) \right| \leq \kappa''.
\end{equation*}
The random variables
\begin{equation*}
Z_{n, \, i}(\omega) = (p_n+1) \frac{h_n^d}{\displaystyle\sup_{B(\omega, \, h_n)} g^{p_n}} \, S_{n, \, i}(\omega), \ i=1, \, \ldots, \, n
\end{equation*}
are therefore uniformly bounded, centered, independent and identically distributed. Let
\begin{eqnarray*}
\tau_n(\omega) &:=& \frac{\varepsilon}{2\kappa''} \, \frac{p_n}{\delta_n \, w_n} \, \frac{n\, \mu_{p_n+1}(\omega) \, h_n^d}{\displaystyle\sup_{B(\omega, \, h_n)} g^{p_n}} \\ 
\textrm{and} \quad \lambda_n(\omega) &:=& \frac{\varepsilon \kappa''}{2} \, \frac{p_n}{\delta_n \, w_n} \, \frac{\mu_{p_n+1}(\omega) \, h_n^d}{\displaystyle\sup_{B(\omega, \, h_n)} g^{p_n}} \, \frac{1}{\Var(Z_{n, \, 1}(\omega))} \\
								&=& \frac{\varepsilon \kappa''}{2} \, \frac{p_n}{\delta_n \, w_n} \, \mu_{p_n+1}(\omega) \, \sup_{B(\omega, \, h_n)} g^{p_n} \, \frac{(p_n+1)^{-2} \, h_n^{-d}}{\EE|S_{n, \, 1}(\omega)|^2}.
\end{eqnarray*}
Recalling~(\ref{defzeta1}), Bernstein's inequality yields, for all $\varepsilon>0$ and $n$ large enough, 
\begin{equation*}
\forall \, \omega\in\Omega_n, \ \PP\left( \delta_n \, \frac{w_n}{p_n} \left| \zeta_n^{(1)}(\omega) \right| > \frac{\varepsilon}{2} \right) \leq 2\exp\left( -\frac{\tau_n(\omega) \lambda_n(\omega)}{2(1+\lambda_n(\omega)/3)} \right).
\end{equation*}
Proposition~\ref{propequivmu}, equation~(\ref{a2}) and the equality $w_n=p_n\, v_n$ entail
\begin{equation*}
\inf_{\omega\in \Omega_n} \tau_n(\omega) \geq \frac{\varepsilon}{4\kappa''} \, \inf_{\Omega} g \, \frac{\sqrt{\log n}}{\delta_n} \sqrt{n \, p_n^{-\overline{\alpha}} \, h_n^d}
\end{equation*}
for large enough $n$. Moreover, straightforward computations yield 
\begin{equation*}
\forall \, y\in [0, \, 1], \ \sup_{x\in \Omega} \, \sup_{u \in B} \left| -1 + \frac{\mu_{p_n}(x)}{\mu_{p_n+1}(x)} \, g(x-h_n u) \, y \right| \leq (1-y) + \frac{\alpha(x) y + \nu_n(y)}{p_n},
\end{equation*}
with $\nu_n$ being a sequence of Borel functions converging uniformly to 0.
Lemma~\ref{lemsq} thus shows that
\begin{equation*}
\sup_{x\in\Omega} \left| \frac{\EE|S_{n, \, 1}(x)|^2}{g^{2 p_n}(x) \, p_n^{-\alpha(x)-2} \, h_n^{-d}} \right|=\operatorname{O}(1) \quad \textrm{as } n\to\infty.
\end{equation*}
Consequently, applying Proposition~\ref{propequivmu} to $\mu_{p_n+1}(\omega)$ entails
\begin{equation*}
\sup_{\omega\in \Omega_n} \frac{1}{\lambda_n(\omega)} = \operatorname{O}\left( \delta_n \sqrt{\frac{n \, p_n^{-\overline{\alpha}} \, h_n^d}{\log n}} \right)
\end{equation*}
as $n\to\infty$. Thus, using once again the fact that the function $x\mapsto 1/[2(x+1/3)]$ is decreasing on $\RR_+$, we get that there exists a constant $c_{\varepsilon}>0$ such that for all $n$ large enough, 
\begin{equation*}
\forall \, \omega\in\Omega_n, \ \PP\left( \delta_n \, \frac{w_n}{p_n} \left| \zeta_n^{(1)}(\omega) \right| > \frac{\varepsilon}{2} \right) \leq 2\exp\left( -c_{\varepsilon} \, \frac{\log n}{\delta_n^2} \right).
\end{equation*}
As a consequence, $\sum_n T_{2, \, n}$ converges and (\ref{random}) is proven: applying Lemma~\ref{lemtech2} completes the proof of Theorem~\ref{urateconsist}. \CQFD
%
%-----------------------------------------------------------------------------
\section*{Appendix: Auxiliary results and proofs}
%-----------------------------------------------------------------------------
%
\noindent
The first lemma of this section is a topological result which shall be needed in several proofs. 
\begin{Lem}
\label{lemtop}
There exists $\beta>0$ such that for every $x\in\Omega$, $B(x, \, \beta)\subset E$.
\end{Lem}
\noindent
{\bf Proof of Lemma~\ref{lemtop}.} Let $U$ denote the interior of $E$ and $\partial E=E\setminus U$ be the (topological) boundary of $E$. Note that $\partial E$ is a closed set since it is the intersection of two closed sets in $\RR^d$; since $\Omega$ is a compact set and $\partial E$ is a closed  set with $\Omega\cap \partial E=\varnothing$, it holds that 
\begin{equation}
\label{dist}
\exists \, \beta>0, \ d(\Omega, \, \partial E) :=\inf_{x\in \Omega} \, \inf_{e \in \partial E} \| x-e \| = 2\beta>0. 
\end{equation}
We shall now prove the result. Pick $x\in \Omega$. If one could find $y\in B(x, \, \beta) \cap E^c$ -- where $E^c$ is the complement of the set $E$ -- then the real number
$$
t_0=\inf\{ t\in [0, \, 1] \, | \, z_t := (1-t)x+ty \notin E \}
$$
would belong to $(0, \, 1)$ since $x\in U$ and $y\in E^c$ which are both open sets. Therefore, because for every $t\in (0, \, t_0)$, $z_t\in E$ and there exists a nonincreasing sequence $(t_k)$ converging to $t_0$ such that $(z_{t_k})\subset E^c \subset U^c$ which is a closed set, one has 
$$
z_{t_0} = \lim_{t \uparrow t_0} z_t \in E \ \mbox{ and } \ z_{t_0} = \lim_{k\to\infty} z_{t_k} \in U^c.
$$
Hence $z_{t_0}\in \partial E$, but $\| x-z_{t_0} \| = t_0 \| x-y \| < \beta$, which contradicts~(\ref{dist}): Lemma~\ref{lemtop} is proven. \CQFD

\vskip2ex
\noindent
We proceed with a technical result we shall need to examine the properties of $m_{p_n}(x)$ and $\mu_{p_n}(x)$ in Lemma~\ref{lemuratiodet} below. It essentially shows that the computation of a conditional high order moment is controlled by the behavior of the conditional survival function $\overline{\mathcal{F}}(\cdot \, | \, x)$ in a neighborhood of 1. 
\begin{Lem}
\label{lemtech1}
Let $h$ be a positive bounded Borel function on $(0, \, 1)$, and let $p_n\to\infty$. If $(NP_1-NP_2)$ hold, then for all $\varepsilon\in (0, \, 1-y_0)$,
\begin{equation*}
 \sup_{x\in \Omega} \left| \frac{\displaystyle\int_{1-\varepsilon}^1 y^{p_n-1} \, h(y) \, \overline{\mathcal{F}}( y \, | \, x) \, dy }{\displaystyle\int_0^1 y^{p_n-1} \, h(y) \, \overline{\mathcal{F}}( y \, | \, x) \, dy} - 1 \right| \to 0 \quad \textrm{as } n\to\infty.
\end{equation*}
\end{Lem}
\noindent
{\bf Proof of Lemma~\ref{lemtech1}.} Let $\varepsilon\in (0, \, 1-y_0)$, $x\in \Omega$ and
consider the expansion
\begin{equation*}
\int_0^1 y^{p_n-1} \, h(y) \, \overline{\mathcal{F}}(y \, | \, x) \, dy = \int_{1-\varepsilon}^1 y^{p_n-1} \, h(y) \, \overline{\mathcal{F}}(y \, | \, x) \, dy \left[1+\frac{ \displaystyle\int_0^{1-\varepsilon} y^{p_n-1} \, h(y) \, \overline{\mathcal{F}}(y \, | \, x) \, dy }{\displaystyle\int_{1-\varepsilon}^1 y^{p_n-1} \, h(y) \, \overline{\mathcal{F}}(y \, | \, x) \, dy } \right].
\end{equation*}
Since, for all $y\in [1-\varepsilon, \, 1]$, the function $x\mapsto \overline{\mathcal{F}}(y \, | \, x)$ is positive and continuous on $\Omega$, it is clear that $\displaystyle\inf_{x\in \Omega} \overline{\mathcal{F}}(y \, | \, x) >0$. Consequently
\begin{eqnarray*}
0\leq \sup_{x\in \Omega} \left| \frac{ \displaystyle\int_0^{1-\varepsilon} y^{p_n-1} \, h(y) \, \overline{\mathcal{F}}(y \, | \, x) \, dy }{\displaystyle\int_{1-\varepsilon}^1 y^{p_n-1} \, h(y) \, \overline{\mathcal{F}}(y \, | \, x) \, dy } \right| & \leq & \sup_{x\in \Omega} \left| \frac{ (1-\varepsilon) \, \displaystyle\sup_{(0, \, 1)} h} {\displaystyle\int_{1-\varepsilon}^1 \left[\frac{y}{1-\varepsilon}\right]^{p_n-1} \, h(y) \, \overline{\mathcal{F}}(y \, | \, x) \, dy } \right| \\ 
 & \leq & \frac{ (1-\varepsilon) \, \displaystyle\sup_{(0, \, 1)} h }{\displaystyle\left[\frac{1-\varepsilon/2}{1-\varepsilon}\right]^{p_n-1} \int_{1-\varepsilon/2}^1 h(y) \, \displaystyle\inf_{x\in \Omega} \overline{\mathcal{F}}(y \, | \, x) \, dy}.
\end{eqnarray*}
Remarking that $\left[\dfrac{1-\varepsilon/2}{1-\varepsilon}\right]^{p_n-1} \to \infty$ as $n\to\infty$, we get the desired result. \CQFD

\vskip2ex
\noindent
The following lemma examines the behavior of the conditional high order moment $m_{p_n}(x)$ and its smoothed version $\mu_{p_n}(x)$ in the nonparametric context.
\begin{Lem}
\label{lemuratiodet}
Assume that $(NP_1-NP_3)$ and $(A_1)$ hold. Let $K$ be a probability density function on $\RR^d$ with support included in $B$. If $p_n\to\infty$ and $p_n \, h_n^{\eta_g} \to 0$ as $n\to \infty$, then
\begin{enumerate}[label=(\roman{*})]
\item $\displaystyle\sup_{x\in \Omega}\left| \frac{\mu_{p_n}(x)}{f(x) \, m_{p_n}(x)} -1 \right| \to 0 \quad \textrm{as } n\to\infty$,
\item $\displaystyle\sup_{x\in \Omega}\left| \frac{m_{p_n+1}(x)}{m_{p_n}(x)} -g(x) \right| \to 0 \quad \textrm{as } n\to\infty$,
\item $\displaystyle\sup_{x\in \Omega}\left| \frac{\mu_{p_n+1}(x)}{\mu_{p_n}(x)} -g(x) \right| \to 0 \quad \textrm{as } n\to\infty$.
\end{enumerate}
\end{Lem}
\noindent
{\bf Proof of Lemma~\ref{lemuratiodet}.} Before starting the proof of this result, use Lemma~\ref{lemtop} to pick $n$ large enough such that $B(x, \, h_n)\subset E$ for all $x\in\Omega$.

\vskip1ex
\noindent
(i) Let us remark that
\begin{equation*}
\mu_{p_n}(x)=\int_{B} K(u) \, f(x-h_n u) \, m_{p_n}(x-h_n u) \, du, 
\end{equation*}
so that 
\begin{equation*}
\frac{\mu_{p_n}(x)}{f(x) \, m_{p_n}(x)}=\int_{B} K(u) \, \frac{f(x-h_n u)}{f(x)} \, \frac{m_{p_n}(x-h_n u)}{m_{p_n}(x)} \, du.
\end{equation*}
Besides,
\begin{equation*}
\frac{m_{p_n}(x-h_n u)}{m_{p_n}(x)} = \frac{g^{p_n}(x-h_n u)}{g^{p_n}(x)} \, \frac{ \displaystyle\int_0^{1} y^{p_n-1} \overline{\mathcal{F}}(y \, | \, x-h_n u) \, dy }{\displaystyle\int_0^1 y^{p_n-1} \overline{\mathcal{F}}(y \, | \, x) \, dy }.
\end{equation*}
From (\ref{a1}), (\ref{a2}) and hypothesis $(NP_3)$, it follows that
\begin{equation*}
\frac{\mu_{p_n}(x)}{f(x) \, m_{p_n}(x)} \to \int_B K(u) \, du = 1
\end{equation*}
uniformly in $x\in \Omega$ as $n\to\infty$, which proves (i). 

\vskip1ex
\noindent
(ii) Similarly, we have 
\begin{equation*}
\frac{m_{p_n+1}(x)}{m_{p_n}(x)} = g(x) \left[ 1+\frac{1}{p_n} \right] \frac{ \displaystyle\int_0^1 y^{p_n} \overline{\mathcal{F}}(y \, | \, x) \, dy }{\displaystyle\int_0^1 y^{p_n-1} \overline{\mathcal{F}}(y \, | \, x) \, dy }.
\end{equation*}
Note that 
\begin{equation*}
1-\frac{ \displaystyle\int_0^1 y^{p_n} \overline{\mathcal{F}}(y \, | \, x) \, dy }{\displaystyle\int_0^1 y^{p_n-1} \overline{\mathcal{F}}(y \, | \, x) \, dy } = \frac{ \displaystyle\int_0^1 y^{p_n-1} (1-y) \overline{\mathcal{F}}(y \, | \, x) \, dy }{\displaystyle\int_0^1 y^{p_n-1} \overline{\mathcal{F}}(y \, | \, x) \, dy }
\end{equation*}
and let $\varepsilon\in (0, \, 1-y_0)$. Lemma~\ref{lemtech1} shows that, for all $n$ large enough,
\begin{equation*}
\sup_{x\in \Omega} \left| 1-\frac{ \displaystyle\int_0^1 y^{p_n} \overline{\mathcal{F}}(y \, | \, x) \, dy }{\displaystyle\int_0^1 y^{p_n-1} \overline{\mathcal{F}}(y \, | \, x) \, dy } \right| \leq (1+\varepsilon) \, \sup_{x\in\Omega} \left| \frac{ \displaystyle\int_{1-\varepsilon}^1 y^{p_n-1} (1-y) \overline{\mathcal{F}}(y \, | \, x) \, dy }{\displaystyle\int_{1-\varepsilon}^1 y^{p_n-1} \overline{\mathcal{F}}(y \, | \, x) \, dy } \right| \leq \varepsilon(1+\varepsilon)
\end{equation*}
and the result follows.

\vskip1ex
\noindent
(iii) is a consequence of (i) and (ii). \CQFD

\vskip2ex
\noindent
The fourth lemma of this section establishes a uniform control of the relative oscillation of $\mu_{p_n}$.
\begin{Lem}
\label{lemosc}
Assume that $(NP_1-NP_3)$, $(K)$ and $(A_1)$ hold. Let $(\varepsilon_n)$ be a sequence of positive real numbers such that $\varepsilon_n \leq h_n$. If $p_n\to\infty$ and $p_n \, h_n^{\eta_g} \to 0$ as $n\to \infty$, then
\begin{equation*}
\sup_{x\in \Omega} \, \sup_{z\in B(x, \, \varepsilon_n)} \left| \frac{\mu_{p_n}(z)}{\mu_{p_n}(x)}-1 \right| = \operatorname{O}\left( \left[ \frac{\varepsilon_n}{h_n} \right]^{\eta_K} \right).
\end{equation*}
\end{Lem}
\noindent
{\bf Proof of Lemma~\ref{lemosc}.} For all $x\in \Omega$ and $z\in B(x, \, \varepsilon_n)$, we have
\begin{equation*}
\left| \mu_{p_n}(x)-\mu_{p_n}(z) \right| \leq \EE\left( Y^{p_n} \left| K_{h_n}(x-X) - K_{h_n}(z-X) \right| \right).
\end{equation*}
Hypothesis $(K)$ and the inclusion $B(z, \, h_n) \subset B(x, \, 2h_n)$ now entail
\begin{eqnarray*}
\left| K_{h_n}(x-X) - K_{h_n}(z-X) \right| & \leq & \frac{c_K}{h_n^d} \, \left[ \frac{\| x-z \|}{h_n} \right]^{\eta_K} \ind_{\{ X\in B(x, \, h_n) \cup B(y, \, h_n) \}} \\ 
										   & \leq & \frac{c_K}{h_n^d} \, \left[ \frac{\varepsilon_n}{h_n} \right]^{\eta_K} \ind_{\{ X\in B(x, \, 2h_n) \}}.
\end{eqnarray*}
Let $\mathcal{V}$ be the volume of the unit ball in $\RR^d$, $\mathcal{K}=\ind_{B} / \mathcal{V} $ be the uniform kernel on $\RR^d$ and let $\mathcal{K}_h(u)=h^{-d} \, \mathcal{K}(u/h)$. The oscillation of $\mu_{p_n}(x)$ is controlled as
\begin{equation}
\label{oscmu}
\sup_{z\in B(x, \, \varepsilon_n)} \left| \mu_{p_n}(x)-\mu_{p_n}(z) \right| \leq 2^d \, c_K \mathcal{V} \ \EE\left( Y^{p_n} \, \mathcal{K}_{2h_n}(x-X) \right) \left[ \frac{\varepsilon_n}{h_n} \right]^{\eta_K}.
\end{equation}
Note that $\mathcal{K}$ is a probability density function on $\RR^d$ with support included in $B$. Therefore, Lemma~\ref{lemuratiodet}(i) yields
\begin{equation*}
\sup_{x\in \Omega}\left| \frac{\EE\left( Y^{p_n} \, \mathcal{K}_{2h_n}(x-X) \right)}{f(x) \, m_{p_n}(x)} -1 \right| \to 0 \quad \textrm{as } n\to\infty.
\end{equation*}
Applying Lemma~\ref{lemuratiodet}(i) once again gives
\begin{equation*}
\sup_{x\in \Omega}\left| \frac{\EE\left( Y^{p_n} \, \mathcal{K}_{2h_n}(x-X) \right)}{\mu_{p_n}(x)} -1 \right| \to 0 \quad \textrm{as } n\to\infty
\end{equation*}
which, together with~(\ref{oscmu}), yields the result. \CQFD

\vskip2ex
\noindent
Lemma~\ref{lemstirunif} below is a useful tool in establishing uniform expansions for ratios of Gamma functions:
\begin{Lem}
\label{lemstirunif}
For all $z, \, z'>0$, one has
\begin{equation*}
\log \frac{\Gamma(z)}{\Gamma(z')} =\left( z-\frac{1}{2} \right) \log z -  \left( z'-\frac{1}{2} \right) \log z' - (z-z') + \operatorname{O}\left( \left| \frac{1}{z}-\frac{1}{z'} \right| \right).
\end{equation*}
\end{Lem}
\noindent
{\bf Proof of Lemma~\ref{lemstirunif}.} From (6.1.50) in Abramovitz and Stegun (1965), p.258, one has
\begin{equation*}
\log \Gamma(z) = \left( z-\frac{1}{2} \right) \log z - z +\frac{1}{2} \log 2\pi + 2\int_0^{\infty} \frac{\arctan(t/z)}{e^{2\pi t}-1} \, dt.
\end{equation*}
Now, since $x\mapsto \arctan x$ is a Lipschitz function on $\RR$, it follows that
\begin{equation*}
\left| \int_0^{\infty} \frac{\arctan(t/z)}{e^{2\pi t}-1} \, dt - \int_0^{\infty} \frac{\arctan(t/z')}{e^{2\pi t}-1} \, dt \right| \leq \left| \frac{1}{z} - \frac{1}{z'} \right| \int_0^{\infty} \frac{t}{e^{2\pi t}-1} \, dt 
\end{equation*}
Remarking that the integral on the right-hand side is convergent yields
\begin{equation*}
\left| \int_0^{\infty} \frac{\arctan(t/z)}{e^{2\pi t}-1} \, dt - \int_0^{\infty} \frac{\arctan(t/z')}{e^{2\pi t}-1} \, dt \right| =\operatorname{O}\left( \left| \frac{1}{z}-\frac{1}{z'} \right| \right)
\end{equation*}
and the result follows. \CQFD

\vskip2ex
\noindent
The next result of this section is a generalisation of Lemma~2 in Girard {\it et al.} (2013). It provides a uniform expansion of $M_n(p_n, \, x)$, see~(\ref{notation_M}) in the proof of Proposition~\ref{propequivmu}, which is the key to the proof of Proposition~\ref{propdeux}.
\begin{Lem}
\label{lemdvpt}
Assume that $(K)$ and $(A_1-A_2)$ hold. For all $x\in \Omega$, $u \in B$ and $n \in \NN\setminus\{ 0 \}$, let 
\begin{eqnarray*} 
\mathcal{L}_n(p_n, \, x, \, u) \! \! &=& \! \! \dfrac{(fC)(x-h_n u) \, \Gamma(\alpha(x-h_n u)+1)}{(fC)(x) \, \Gamma(\alpha(x)+1)} \exp\left[ p_n \frac{\Delta_n^g(x, \, u)}{g(x)} - \log(p_n) \Delta_n^{\alpha}(x, \, u) \right], \\
\Lambda_n(p_n, \, x) \! \! &=& \! \! \dfrac{M_n(p_n, \, x)}{ f(x) \, C(x) \, g^{p_n}(x) }.
\end{eqnarray*}
If $p_n\to\infty$ and $p_n \, h_n^{\eta_g} \to 0$, then
\begin{equation*}
\sup_{x\in\Omega} \left| \frac{\Lambda_n(p_n, \, x)}{\alpha(x) b(p_n+1, \, \alpha(x)) } - 1 \right| \to 0
\end{equation*}
and 
\begin{equation*}
\sup_{x\in\Omega} \left| \frac{\Lambda_n(p_n, \, x)}{\alpha(x) b(p_n+1, \, \alpha(x)) } - \int_B \mathcal{L}_n(p_n, \, x, \, u) \,  K(u) \, du \right| = \operatorname{O}\left( h_n^{\eta_g} \lor \frac{h_n^{\eta_{\alpha}}}{p_n} \right).
\end{equation*}
\end{Lem}
\noindent
{\bf Proof of Lemma~\ref{lemdvpt}.} Using Lemma~\ref{lemtop}, we can pick $n$ large enough such that $B(x, \, h_n)\subset E$ for all $x\in\Omega$. Introducing
\begin{equation}
\label{defqn}
Q_n(x, \, u)=\frac{(fC)(x-h_n u) \, \Gamma(\alpha(x-h_n u)+1)}{(fC)(x) \, \Gamma(\alpha(x)+1)},
\end{equation}
we have
\begin{equation}
\label{ref_lm}
\dfrac{\Lambda_n(p_n, \, x)}{\alpha(x) \, b(p_n+1, \, \alpha(x))} = \int_{B} Q_n(x, \, u) \, \frac{\Gamma(p_n+1+\alpha(x))}{\Gamma( p_n+1+\alpha(x-h_n u))} \, \frac{g^{p_n}(x-h_n u)}{g^{p_n}(x)} \, K(u) \, du.
\end{equation}
The set $\Omega$ being a compact set, the set $\Omega'=\{x'\in \RR^d \, | \, \exists \, x\in \Omega, \, \| x-x' \| \leq h_n \}$ is compact as well, and $\Omega'\subset E$: since $f, \, C$ and $\alpha$ are continuous on the compact set $\Omega'\subset E$, they are uniformly continuous on $\Omega'$. Furthermore, since $\alpha$ is bounded on $\Omega'$ and $\Gamma$ is continuous on $(0, \, \infty)$, the function $x\mapsto \Gamma(\alpha(x)+1)$ is uniformly continuous on $\Omega'$, so that
\begin{equation}
\label{rapqn}
\sup_{x\in\Omega} \, \sup_{u\in B} \left| Q_n(x, \, u) - 1 \right| \to 0
\end{equation}
as $n\to\infty$. Moreover, since $p_n \, h_n^{\eta_g} \to 0$, we get
\begin{equation*}
\sup_{x\in\Omega} \, \sup_{u \in B} |\log(p_n) \, \Delta_n^{\alpha}(x, \, u)| = \operatorname{O}\left( h_n^{\eta_{\alpha}} |\log p_n| \right) = \operatorname{O}\left( \Big[ h_n^{\eta_g} p_n \Big]^{\eta_{\alpha}/\eta_g} \, \frac{|\log p_n|}{p_n^{\eta_{\alpha}/\eta_g}} \right) \to 0
\end{equation*}
as $n\to\infty$ and Lemma~\ref{lemstirunif} yields
\begin{equation}
\label{rapgam}
\sup_{x\in \Omega} \, \sup_{u\in B} \left| \exp(\log(p_n) \Delta_n^{\alpha}(x, \, u)) \, \frac{\Gamma(p_n+1+\alpha(x))}{\Gamma(p_n+1+\alpha(x-h_n u))}-1 \right| = \operatorname{O}\left( \frac{h_n^{\eta_{\alpha}}}{p_n} \right).
\end{equation}
Besides,
\begin{equation}
\label{rapg}
\frac{g^{p_n}(x-h_n u)}{g^{p_n}(x)} = \exp\left[ p_n \log\left(1+\frac{\Delta_n^g(x, \, u)}{g(x)} \right) \right] 
\end{equation}
where
\begin{equation*}
\sup_{x\in\Omega} \, \sup_{u \in B} p_n \, \left| \frac{\Delta_n^g(x, \, u)}{g(x)} \right| \to 0
\end{equation*}
as $n\to\infty$, see (\ref{myeq}). Replacing (\ref{rapqn}), (\ref{rapgam}) and (\ref{rapg}) in (\ref{ref_lm}) gives both results. \CQFD 

\vskip2ex
\noindent
The aim of Lemma~\ref{lemlin} below is to linearise the random variable $\xi_n(x)$ appearing in the proof of Theorem~\ref{urateconsist}:  

\begin{Lem}
\label{lemlin} 
The random variable $\xi_n(x)$ can be expanded as
\begin{equation*}
\xi_n(x) = \frac{1}{a p_n} \left[ \zeta_n^{(1)}(x) - \zeta_n^{(2)}(x) + \left(\frac{\mu_{p_n+1}(x)}{\widehat{\mu}_{p_n+1}(x)}-1 \right) \zeta_n^{(1)}(x) - \left( \frac{\mu_{(a+1) p_n+1}(x)}{\widehat{\mu}_{(a+1) p_n+1}(x)}-1 \right) \zeta_n^{(2)}(x) \right]
\end{equation*}
where 
\begin{eqnarray*}
\zeta_n^{(1)}(x) \! \! &=& \! \! (p_n+1) \frac{\mu_{p_n}(x)}{\mu_{p_n+1}(x)} \left[ \frac{\widehat{\mu}_{p_n+1}(x)}{\mu_{p_n+1}(x)} - \frac{\widehat{\mu}_{p_n}(x)}{\mu_{p_n}(x)} \right] \\
\textrm{and } \ \zeta_n^{(2)}(x) \! \! &=& \! \! [(a+1)p_n+1] \frac{\mu_{(a+1)p_n}(x)}{\mu_{(a+1)p_n+1}(x)} \left[ \frac{\widehat{\mu}_{(a+1)p_n+1}(x)}{\mu_{(a+1)p_n+1}(x)} - \frac{\widehat{\mu}_{(a+1)p_n}(x)}{\mu_{(a+1)p_n}(x)} \right].
\end{eqnarray*}
\end{Lem}
\noindent
{\bf Proof of Lemma~\ref{lemlin}.} Straightforward computations yield
\begin{equation}
\label{debutlemlin}
ap_n \, \xi_n(x) = D_n^{(1)}(x)-D_n^{(2)}(x) 
\end{equation}
with
\begin{eqnarray*}
D_n^{(1)}(x) &:=& (p_n+1) \frac{\mu_{p_n}(x)}{\mu_{p_n+1}(x)} \, \frac{\mu_{p_n+1}(x)}{\widehat{\mu}_{p_n+1}(x)} \, \left[ \frac{\widehat{\mu}_{p_n+1}(x)}{\mu_{p_n+1}(x)} - \frac{\widehat{\mu}_{p_n}(x)}{\mu_{p_n}(x)}  \right], \\
D_n^{(2)}(x) &:=& [(a+1)p_n+1] \frac{\mu_{(a+1)p_n}(x)}{\mu_{(a+1) p_n+1}(x)} \, \frac{\mu_{(a+1) p_n+1}(x)}{\widehat{\mu}_{(a+1) p_n+1}(x)} \, \left[ \frac{\widehat{\mu}_{(a+1)p_n+1}(x)}{\mu_{(a+1)p_n+1}(x)} - \frac{\widehat{\mu}_{(a+1)p_n}(x)}{\mu_{(a+1)p_n}(x)} \right].
\end{eqnarray*}
This leads to
\begin{eqnarray*}
D_n^{(1)}(x) = \frac{\mu_{p_n+1}(x)}{\widehat{\mu}_{p_n+1}(x)} \, \zeta_n^{(1)}(x) \quad \textrm{and} \quad D_n^{(2)}(x) = \frac{\mu_{(a+1) p_n+1}(x)}{\widehat{\mu}_{(a+1) p_n+1}(x)} \, \zeta_n^{(2)}(x);
\end{eqnarray*}
replacing in (\ref{debutlemlin}) concludes the proof of Lemma~\ref{lemlin}. \CQFD

\vskip2ex
\noindent
We shall next take a closer look at the behavior of the functions $\varepsilon_n(p_n, \, x)$, see~(\ref{devmu}) in the proof of Proposition~\ref{propequivmu}. We first introduce some tools necessary for this study. For an arbitrary set $\mathcal{S}$, $\mathcal{F}(\mathcal{S})$ is the set of all sequences of functions $u_n: \NN \times \mathcal{S} \to \RR$, denoted by $u_n(t, \, x)$. Let $\mathcal{C}(\mathcal{S}) \subset \mathcal{F}(\mathcal{S})$ be the subset of all the elements $u\in \mathcal{F}(\mathcal{S})$ such that $u$ meets the following requirements: 

\vskip1ex

{\bf  $(Q_1)$} There exists $N_1\in \NN$ such that for all $t\in \NN$, $\displaystyle\sup_{n\geq N_1} \, \sup_{x\in \mathcal{S}} |u_n(t, \, x)| < \infty$.

\vskip1ex

{\bf  $(Q_2)$} There exists $N_2\in \NN$ such that for all $t, \, t'\in \NN$, $p_n \, \displaystyle\sup_{n\geq N_2} \, \sup_{x\in \mathcal{S}} | u_n(t', \, x)-u_n(t, \, x) | < \infty$.

\vskip1ex
\noindent 
Finally, $\mathcal{D}(\mathcal{S})$ is a subset of $\mathcal{C}(\mathcal{S})$ whose elements are bounded from below:
\begin{equation*}
\mathcal{D}(\mathcal{S}) = \{ u \in \mathcal{C}(\mathcal{S}) \, | \, \exists \, N_0 \in \NN, \ \forall \, t\in \NN, \ \exists \, M(t)>0, \ \inf_{n\geq N_0} \, \inf_{x\in \mathcal{S}} u_n(t, \, x) \geq M(t) \}.
\end{equation*}
Lemma~\ref{lemalg} lists some properties of the sets $\mathcal{C}(\mathcal{S})$ and $\mathcal{D}(\mathcal{S})$.
\begin{Lem}
\label{lemalg}
Let $\mathcal{S}$ be an arbitrary set. Then:
\begin{enumerate}[label=(\roman{*})]
\item $\mathcal{C}(\mathcal{S})$ is a linear subspace of $\mathcal{F}(\mathcal{S})$ which is closed under multiplication.
\item $\mathcal{D}(\mathcal{S})$ is closed under multiplication and division. 
\item Let $u\in \mathcal{F}(\mathcal{S})$ such that there exists a sequence of uniformly bounded real functions $(\delta_n)$ on $\mathcal{S}$ with
\begin{equation*}
\forall \, t\in \NN, \ \sup_{x\in \mathcal{S}} \left| u_n(t, \, x)- \left[ 1+\frac{\delta_n(x)}{p_n+t} \right] \right| = \operatorname{o}\left( \frac{1}{p_n} \right).
\end{equation*}
Then $u \in \mathcal{D}(\mathcal{S})$.
\item If $(\mathcal{S}', \, \mathcal{T}, \, \mu)$ is a finite measure space and if $u\in \mathcal{C}(\mathcal{S} \times \mathcal{S}')$ (resp. $\mathcal{D}(\mathcal{S} \times \mathcal{S}')$) is such that $x'\mapsto u_n(t, \, (x, \, x'))$ is measurable for every $t\in \NN$ and $x\in \mathcal{S}$, then
\begin{equation*}
(n, \, t, \, x) \mapsto \int_{\mathcal{S}'} u_n(t, \, (x, \, x')) \, \mu(dx') \in \mathcal{C}(\mathcal{S}) \quad \textrm{(resp. $\mathcal{D}(\mathcal{S})$)}.
\end{equation*}
\end{enumerate}
\end{Lem}
\noindent
{\bf Proof of Lemma~\ref{lemalg}.} (i) Since it is straightforward that $\mathcal{C}(\mathcal{S})$ is a linear subspace of $\mathcal{F}(\mathcal{S})$, it is enough to prove that $\mathcal{C}(\mathcal{S})$ is closed under multiplication. Let $u, \, v \in \mathcal{C}(\mathcal{S})$ and let $w_n(t, \, x)=u_n(t, \, x) \, v_n(t, \, x)$. One has, for all $x\in \mathcal{S}$ and $t, \, t'\in \NN$:
\begin{equation*}
w_n(t', \, x)-w_n(t, \, x) = u_n(t', \, x)[v_n(t', \, x)-v_n(t, \, x)] +v_n(t, \, x)[u_n(t', \, x)-u_n(t, \, x)].
\end{equation*}
Since $u$ and $v$ satisfy requirements $(Q_1)$ and $(Q_2)$, this equality therefore shows that $w$ satisfies $(Q_2)$, and (i) is proven. 

\vskip1ex
\noindent
(ii) Stability under multiplication is a direct consequence of (i). It is then enough to prove that if $u\in \mathcal{D}(\mathcal{S})$, then $1/u\in \mathcal{D}(\mathcal{S})$. Let $w=1/u$: $w$ clearly satisfies $(Q_1)$ and for all $t\in \NN$ and $n$ large enough, $\displaystyle\inf_{x\in \mathcal{S}} w_n(t, \, x)$ is bounded from below by a positive constant. Finally, for all $t, \, t'\in \NN$, there exists $N_2\in \NN$ such that:
\begin{equation*}
p_n \sup_{n\geq N_2} \, \sup_{x\in \mathcal{S}} \left| \frac{1}{u_n(t, \, x)} - \frac{1}{u_n(t', \, x)} \right| \leq \frac{1}{M(t) \, M(t')} \, p_n \, \sup_{n\geq N_2} \, \sup_{x\in \mathcal{S}} |u_n(t', \, x)-u_n(t, \, x)| <\infty.
\end{equation*}
This is enough to conclude that $w\in \mathcal{C}(\mathcal{S})$, and thus $w\in \mathcal{D}(\mathcal{S})$, which concludes the proof of (ii).

\vskip1ex
\noindent
(iii) Just note that $1/(p_n+t)=1/p_n +\operatorname{o}(1/p_n)$, from which (iii) readily follows. 

\vskip1ex
\noindent
(iv) Let $u\in \mathcal{C}(\mathcal{S} \times \mathcal{S}')$ and
\begin{equation*}
v:(n, \, t, \, x) \mapsto \int_{\mathcal{S}'} u_n(t, \, (x, \, x')) \, \mu(dx') \in \mathcal{F}(\mathcal{S}).
\end{equation*}
Then, for all $t\in \NN$, since $\mu$ is a finite measure on $\mathcal{S}'$, it follows that there exists $N_1\in \NN$ such that
\begin{equation*}
\sup_{n\geq N_1} \, \sup_{x\in \mathcal{S}} | v_n(t, \, x) | \leq \sup_{n\geq N_1} \, \sup_{(x, \, x')\in \mathcal{S}\times \mathcal{S}'} |u_n(t, \, (x, \, x'))| \int_{\mathcal{S}'} \mu(dx') <\infty.
\end{equation*}
Besides, for all $t'\in \NN$, there exists $N_2\in \NN$ such that
\begin{equation*}
p_n \sup_{n\geq N_2} \, \sup_{x\in \mathcal{S}} | v_n(t', \, x)-v_n(t, \, x) | \leq p_n \sup_{n\geq N_2} \, \sup_{(x, \, x')\in \mathcal{S}\times \mathcal{S}'} |u_n(t', \, (x, \, x'))-u_n(t, \, (x, \, x'))| \int_{\mathcal{S}'} \mu(dx') <\infty
\end{equation*}
so that $v\in \mathcal{C}(\mathcal{S})$. If $u\in \mathcal{D}(\mathcal{S} \times \mathcal{S}')$, then there exists $M(t)>0$ and $N_0\in \NN$ such that
\begin{equation*}
\inf_{n\geq N_0} \, \inf_{x\in \mathcal{S}} v_n(t, \, x) \geq M(t) \int_{\mathcal{S}'} \mu(dx') >0 
\end{equation*}
so that $v\in \mathcal{D}(\mathcal{S})$, and (iv) is proven. \CQFD

\vskip2ex
\noindent
Lemma~\ref{lemerr} below essentially gives the order of magnitude of $M_n(p_n+t, \, x)$ and the error term $E_n(p_n+t, \, x)$ in the expansion of $\mu_{p_n}(x)$:
\begin{Lem}
\label{lemerr}
Assume that $(A_1-A_2)$ hold. If $p_n\to\infty$ and $p_n \, h_n^{\eta_g} \to 0$ as $n\to \infty$, then 
\begin{enumerate}[label=(\roman{*})]
\item $(n, \, t, \, x) \mapsto (p_n+t)^{\alpha(x)} \, \dfrac{M_n(p_n+t, \, x)}{g^{p_n+t}(x)} \in \mathcal{D}(\Omega)$;
\item $(n, \, t, \, x) \mapsto (p_n+t)^{[\alpha+\beta](x)} \, E_n(p_n+t, \, x) \in \mathcal{C}(\Omega)$.
\end{enumerate}
\end{Lem}
\noindent
{\bf Proof of Lemma~\ref{lemerr}.} Before proving this result, note that applying Lemma~\ref{lemtop}, we can pick $n$ large enough such that $B(x, \, h_n)\subset E$ for all $x\in\Omega$. 

\vskip1ex
\noindent
(i) Recalling the notations of Lemma~\ref{lemdvpt}, we have
\begin{equation*}
\frac{M_n(p_n, \, x)}{\alpha(x) \, b(p_n+1, \, \alpha(x)) \, (f C g^{p_n})(x)}  =  \int_{B} Q_n(x, \, u) \, \frac{\Gamma(p_n+1+\alpha(x))}{\Gamma( p_n+1+\alpha(x-h_n u))} \, \frac{g^{p_n}(x-h_n u)}{g^{p_n}(x)} \, K(u) \, du.
\end{equation*}
Since
\begin{equation*}
(p_n+t)^{\alpha(x)} \, b(p_n+t+1, \, \alpha(x)) = \frac{(p_n+t)^{\alpha(x)+1}}{\alpha(x)} b(p_n+t, \, \alpha(x)+1),
\end{equation*}
Lemma~\ref{lemstirunif} and  Lemma~\ref{lemalg}(iii) yield
\begin{equation}
\label{betaE}
(n, \, t, \, x) \mapsto (p_n+t)^{\alpha(x)} \, b(p_n+t+1, \, \alpha(x)) \in \mathcal{D}(\Omega). 
\end{equation}
Consequently, it is enough to show that 
\begin{equation*}
(n, \, t, \, x) \mapsto \int_{B} Q_n(x, \, u) \, \frac{\Gamma(p_n+t+1+\alpha(x))}{\Gamma( p_n+t+1+\alpha(x-h_n u))} \, \frac{g^{p_n+t}(x-h_n u)}{g^{p_n+t}(x)} \, K(u) \, du \in \mathcal{D}(\Omega).
\end{equation*}
From~(\ref{rapgam}) and in view of
\begin{equation*}
\frac{\Gamma(p_n+2+\alpha(x))}{\Gamma(p_n+2+\alpha(x-h_n u))} - \frac{\Gamma(p_n+1+\alpha(x))}{\Gamma(p_n+1+\alpha(x-h_n u))} \! = \! \frac{\Gamma(p_n+1+\alpha(x))}{\Gamma(p_n+1+\alpha(x-h_n u))} \, \frac{-\Delta_n^{\alpha}(x, \, u)}{p_n+1+\alpha(x-h_n u)},
\end{equation*}
it follows by induction that
\begin{equation}
\label{gammaE}
(n, \, t, \, (x, \, u)) \mapsto \frac{\Gamma(p_n+t+1+\alpha(x))}{\Gamma(p_n+t+1+\alpha(x-h_n u))} \in \mathcal{D}(\Omega\times B).
\end{equation}
Then, using the relation
\begin{equation*}
\frac{g^{p_n+1}(x-h_n u)}{g^{p_n+1}(x)} - \frac{g^{p_n}(x-h_n u)}{g^{p_n}(x)} = \frac{g^{p_n}(x-h_n u)}{g^{p_n}(x)} \, \frac{\Delta_n^g(x, \, u)}{g(x)}
\end{equation*}
along with~(\ref{rapg}) gives, by induction, 
\begin{equation}
\label{gE}
(n, \, t, \, (x, \, u)) \mapsto \frac{g^{p_n+t}(x-h_n u)}{g^{p_n+t}(x)} \in \mathcal{D}(\Omega\times B).
\end{equation}
As a consequence of Lemma~\ref{lemalg}iv), (i) is proven. 

\vskip1ex
\noindent
(ii) First and foremost, recall that from~(\ref{notation_E}),
\begin{equation*}
E_n(p_n+t, \, x)=\int_{B} \left( f g^{p_n+t} \right)(x-h_n u) \, (p_n+t) \, \mathcal{I}_{\alpha+\beta, \, D}(p_n+t, \, x-h_n u) \, K(u) \, du.
\end{equation*}
In view of Lemma~\ref{lemalg}(iv), it is then enough to show that 
\begin{equation*}
(n, \, t, \, (x, \, u)) \mapsto (p_n+t)^{[\alpha+\beta](x)+1} \, \frac{g^{p_n+t}(x-h_n u)}{g^{p_n+t}(x)} \, \mathcal{I}_{\alpha+\beta, \, D}(p_n+t, \, x-h_n u) \in \mathcal{C}(\Omega\times B).
\end{equation*}
Using~(\ref{gE}), we shall only prove that
\begin{equation*}
(n, \, t, \, (x, \, u)) \mapsto (p_n+t)^{[\alpha+\beta](x)+1} \, \mathcal{I}_{\alpha+\beta, \, D}(p_n+t, \, x-h_n u) \in \mathcal{C}(\Omega\times B).
\end{equation*}
Since
\begin{equation*}
(p_n+t)^{[\alpha+\beta](x)+1}=p_n^{[\alpha+\beta](x)+1} \, (1+t/p_n)^{[\alpha+\beta](x)+1}
\end{equation*}
and since $(n, \, t, \, x)\mapsto (1+t/p_n)^{[\alpha+\beta](x)+1} \in\mathcal{D}(\Omega)$, in view of Lemma~\ref{lemalg}(i) and (ii), it is sufficient to show the latter property for the function defined by
\begin{equation}
\label{defw}
w_n(t, \, (x, \, u)) = p_n^{[\alpha+\beta](x)+1} \, \mathcal{I}_{\alpha+\beta, \, D}(p_n+t, \, x-h_n u).
\end{equation}
For all $t\in \NN\setminus\{ 0 \}$, let $R_t:[1, \, \infty) \to [0, \, \infty)$ be the function defined by 
\begin{equation*}
\forall \, y\geq 1, \ R_t(y) = y\left\{ 1-\left[ 1-\frac{1}{y} \right]^t \right\}. 
\end{equation*}
For all $t\in \NN\setminus\{ 0 \}$, $R_t$ is a bounded Borel function on $[1, \, \infty)$, and one has, for all $t < t'\in \NN$,
\begin{equation}
\label{diffw}
p_n[w_n(t', \, (x, \, u))-w_n(t, \, (x, \, u))] = -p_n^{[\alpha+\beta](x)+2} \, \mathcal{I}_{\alpha+\beta+1, \, D R_{t'-t}}(p_n+t, \, x-h_n u).
\end{equation}
Remark that for all $j, \, t\in \NN$, $(x, \, u)\in\Omega\times B$ and every bounded Borel function $H$ on $\Omega \times [1, \, \infty)$,
\begin{equation*}
\left| \mathcal{I}_{\alpha+\beta+j, \, H}(p_n+t, \, x-h_n u) \right| \leq b(p_n+t, \, [\alpha+\beta](x-h_n u)+j+1) \, \sup_{\Omega\times[1, \, \infty)} |H|.
\end{equation*}
Finally, Lemma~\ref{lemstirunif} shows that 
\begin{equation*}
\sup_{x\in\Omega} \left| \frac{b(p_n+t, \, [\alpha+\beta](x)+j+1)}{\Gamma([\alpha+\beta](x)+j+1) \, p_n^{-[\alpha+\beta](x)-j-1}} - 1 \right| \to 0.
\end{equation*}
The result follows from~(\ref{defw}) and~(\ref{diffw}). \CQFD

\vskip2ex
\noindent
The next result is particularly useful for providing a uniform asymptotic bound of the second-order moments that appear when computing the rate of convergence in the proof of Theorem~\ref{urateconsist}. This result is an analogue of Lemma~4 in Girard {\it et al.} (2013).
\begin{Lem}
\label{lemsq}
Assume that $(SP)$, $(K)$, $(A_1-A_2)$ hold, $p_n\to\infty$ and $p_n \, h_n^{\eta_g} \to 0$ as $n\to \infty$. Let $(b_{n, \, 0})$ and $(b_{n, \, 1})$ be sequences of Borel functions on $\Omega$ such that there exist sequences of Borel functions $(H_{n, \, 0})$ and $(H_{n, \, 1})$, uniformly bounded on $[0, \, 1]$ with
\begin{equation*}
\forall \, y\in [0, \, 1], \ \sup_{x\in \Omega} \, \sup_{u \in B} \left| b_{n, \, 0}(x) + b_{n, \, 1}(x) \, g(x-h_n u) \, y \right| \leq H_{n, \, 0}(y) \, (1-y) + \frac{H_{n, \, 1}(y)}{p_n}.
\end{equation*}
Then, the sequence of random variables
\begin{equation*}
S_n(x)=Y^{p_n} \left[ b_{n, \, 0}(x) + b_{n, \, 1}(x) \, Y \right] K_{h_n}(x-X)
\end{equation*}
is such that
\begin{equation*}
\sup_{x\in\Omega} \left| \frac{\EE|S_n(x)|^2}{g^{2 p_n}(x) \, p_n^{-\alpha(x)-2} \, h_n^{-d}} \right|=\operatorname{O}(1) \quad \textrm{as } n\to\infty.
\end{equation*}
\end{Lem}
\noindent
{\bf Proof of Lemma~\ref{lemsq}.} Using Lemma~\ref{lemtop}, we can pick $n$ large enough such that $B(x, \, h_n)\subset E$ for all $x\in\Omega$. Conditioning on $X$ yields
\begin{eqnarray*}
\EE|S_n(x)|^2 &=& \int_E \EE \left[ Y^{2 p_n} \left. \left| b_{n, \, 0}(x) + b_{n, \, 1}(x) \, Y \right|^2 \, \right| \, X=v \right] K_{h_n}^2(x-v) \, f(v) \, dv \\
			  &=& h_n^{-d} \int_{B} \EE \left[ Y^{2 p_n} \left. \left| b_{n, \, 0}(x) + b_{n, \, 1}(x) \, Y \right|^2 \, \right| \, X=x-h_n u \right] K^2(u) \, f(x-h_n u) \, du.
\end{eqnarray*}
Now, given $X=x-h_n u$, we have $W_n(x, \, u) := {Y}/{g(x-h_n u)}\leq 1$.  Introducing the bounded sequence
\begin{equation*}
c_n := 2 \sup_{\substack{[0, \, 1] \\ n\in \NN}} \left\{ |H_{n, \, 0}|^2, \, |H_{n, \, 1}|^2 \right\} \, \sup_{x\in\Omega} \, \sup_{u\in B} \left| \frac{g^{2 p_n}(x-h_n u)}{g^{2 p_n}(x)} \right|,
\end{equation*}
H\"older's inequality entails, given $\{ X=x-h_n u \}$,
\begin{equation*}
Y^{2 p_n} \left| b_{n, \, 0}(x) + b_{n, \, 1}(x) \, Y \right|^2 \leq c_n \, g^{2 p_n}(x) \, W_n^{2 p_n}(x, \, u) \left[ (1-W_n(x, \, u))^2 + \frac{1}{p_n^2} \right].
\end{equation*}
It is therefore sufficient to prove that, for all $j\in \NN$:
\begin{equation}
\label{domW}
\sup_{x\in\Omega} \, \sup_{u\in B} \left| \frac{\EE\left( \left. W_n^{2 p_n}(x, \, u) (1-W_n(x, \, u))^{j} \, \right| \, X=x-h_n u \right)}{p_n^{-\alpha(x)-j}} \right| = \operatorname{O}\left( 1 \right). 
\end{equation}
Integrating by parts yields
\begin{eqnarray*}
\EE\left( \left. W_n^{2 p_n}(x, \, u) (1-W_n(x, \, u))^{j} \, \right| \, X=x-h_n u \right) &=& \int_0^1 \frac{d}{dy} \left[ y^{2 p_n} \, (1-y)^j \right] \, \overline{\mathcal{F}}(y \, | \, x-h_n u) \, dy \\
																						   &\leq & 2p_n \int_0^1 y^{2 p_n-1} \, (1-y)^j \, \overline{\mathcal{F}}(y \, | \, x-h_n u) \, dy 
\end{eqnarray*}
since, given $\{ X=x-h_n u \}$, $W_n(x, \, u)$ has survival function $\overline{\mathcal{F}}(\cdot \, | \, x-h_n u)$. To conclude, observe that if $\gamma$ is a positive H\"older continuous function on $\RR^d$, then
\begin{equation*}
\int_0^1 y^{2p_n-1} \, \overline{\mathcal{F}}_{\gamma}(y \, | \, x-h_n u) \, dy= b(2p_n, \, \gamma(x-h_n u)+1).
\end{equation*}
From~(\ref{rapgam}) and Stirling's formula, it follows that
\begin{equation*}
\sup_{x\in \Omega} \, \sup_{u\in B} \left| p_n^{\gamma(x)+1} \int_0^1 y^{2p_n-1} \, \overline{\mathcal{F}}_{\gamma}(y \, | \, x-h_n u) \, dy \right|= \operatorname{O}\left( 1 \right)
\end{equation*}
because $p_n \, h_n^{\eta_g} \to 0$ as $n\to\infty$. Finally, for all $y\in [0, \, 1]$,
\begin{equation*}
\overline{\mathcal{F}}(y \, | \, v)=C(v) \, \overline{\mathcal{F}}_{\alpha}(y \, | \, v) + D\left( v, \, (1-y)^{-1} \right) \, \overline{\mathcal{F}}_{\alpha+\beta}(y \, | \, v), 
\end{equation*}
and Lemma \ref{lemerr}(ii) yields~(\ref{domW}), which ends the proof of Lemma~\ref{lemsq}. \CQFD

\vskip2ex
\noindent
The final lemma is the last step in the proof of Theorem~\ref{urateconsist}.
\begin{Lem}
\label{lemtech2}
Let $(X_n)$ be a sequence of positive real-valued random variables such that for every positive nonrandom sequence $(\delta_n)$ converging to 0, the random sequence $(\delta_n X_n)$ converges to 0 almost surely. Then
$$
\PP\left( \limsup_{n\to\infty} X_n = +\infty \right) = 0 \ \mbox{ {\it i.e.} } \ X_n = \operatorname{O}(1) \ \mbox{ almost surely.} 
$$
\end{Lem}
\noindent
{\bf Proof of Lemma~\ref{lemtech2}.} Assume that there exists $\varepsilon>0$ such that $\PP\left( \displaystyle\limsup_{n\to\infty} X_n = +\infty \right) \geq \varepsilon$. Since by definition $\displaystyle\limsup_{n\to\infty} X_n = \lim_{n\to\infty} \, \sup_{p\geq n} X_p$ is the limit of a nonincreasing sequence, one has 
$$
\forall \, k\in \NN, \ \forall \, n\in \NN, \  \PP\left( \bigcup_{p\geq n} \left\{ X_p \geq k \right\} \right) \geq \varepsilon.
$$
From this we deduce
\begin{equation}
\label{constNkNprimek}
\forall \, k\in \NN, \ \forall \, n\in \NN, \ \exists \, n'\geq n, \ \PP\left( \bigcup_{p=n}^{n'} \left\{ X_p \geq k \right\} \right) \geq \varepsilon/2.
\end{equation}
We now build a sequence $(N_k)$ by induction: start by using (\ref{constNkNprimek}) with $k=n=1 =: N_1$ to obtain $N_2 > N_1$ such that
$$
\PP\left( \bigcup_{p=N_1}^{N_2-1} \left\{ X_p \geq 1 \right\} \right) \geq \varepsilon/2.
$$
Then for an arbitrary $k\geq 1$, if $N_k$ is given, apply (\ref{constNkNprimek}) to get $N_{k+1} > N_k$ such that 
$$
\PP\left( \bigcup_{p=N_k}^{N_{k+1}-1} \left\{ X_p \geq k \right\} \right) \geq \varepsilon/2.
$$
The sequence $(N_k)$ is thus an increasing sequence of integers. Let $\delta_n=1/k$ if $N_k \leq n < N_{k+1}$. It is clear that $(\delta_n)$ is a positive sequence which converges to 0. Besides, for all $k \in \NN\setminus\{ 0 \}$ it holds that
$$
\PP\left( \sup_{p\geq N_k} \delta_p X_p \geq 1 \right) = \PP\left( \bigcup_{p\geq N_k} \{ \delta_p X_p \geq 1 \} \right) \geq \PP\left( \bigcup_{p=N_k}^{N_{k+1}-1} \{ \delta_p X_p \geq 1 \} \right) \geq \varepsilon/2.
$$
This entails 
$$
\liminf_{n\to\infty} \PP\left( \sup_{p\geq n} \delta_p X_p \geq 1 \right) \geq \varepsilon/2 >0.
$$
Hence $(\delta_n X_n)$ does not converge almost surely to 0, from which the result follows. \CQFD
\section*{Acknowledgements}
The authors acknowledge the editor and an anonymous referees for their helpful comments which led to significant enhancements of this article.
%
%-----------------------------------------------------------------------------
\section*{References}
\vskip1ex
\noindent
Abramovitz, M., Stegun, I. (1965). {\it Handbook of Mathematical Functions}, Dover. 
\vskip1ex
\noindent
Aragon, Y., Daouia, A., Thomas-Agnan, C. (2005). Nonparametric frontier estimation: a conditional quantile-based approach. {\em Econometric Theory} {\bf 21}(2), 358--389.
\vskip1ex
\noindent
Bingham, N.H., Goldie, C.M., Teugels, J.L. (1987). {\it Regular Variation}, Cambridge, U.K.: Cambridge University Press.
\vskip1ex
\noindent
Cazals, C., Florens, J.-P., Simar, L. (2002). Nonparametric frontier estimation: a robust approach. {\em J. Econometrics} {\bf 106}(1), 1--25.
\vskip1ex
\noindent
Daouia, A., Simar, L. (2005). Robust nonparametric estimators of monotone boundaries. {\it J. Multivariate Anal.} {\bf 96}, 311--331.
\vskip1ex
\noindent 
Deprins, D., Simar, L., Tulkens, H. (1984). Measuring labor efficiency in post offices. In P.~Pestieau, M.~Marchand and H.~Tulkens, editors, {\it The Performance of Public Enterprises: Concepts and Measurements}. North Holland ed, Amsterdam, 243--267.
\vskip1ex
\noindent
Einmahl, U., Mason, D.M. (2000). An empirical process approach to the uniform consistency of kernel-type function estimators. {\it J. Theor. Probab.} {\bf 13}(1), 1--37.
\vskip1ex
\noindent 
Embrechts, P.,  Kl\"uppelberg, C., Mikosch, T. (1997). {\em Modelling extremal events}, Springer.
\vskip1ex
\noindent
Gardes, L., Stupfler, G. (2013). Estimation of the conditional tail-index using a smoothed local Hill estimator. {\it Extremes}, to appear. Available at {\tt http://hal.archives-ouvertes.fr/hal-00739454}.
\vskip1ex
\noindent
Geffroy, J. (1964). Sur un probl\`eme d'estimation g\'eom\'etrique. {\it Publ. Inst. Statist. Univ. Paris} XIII, 191--210.
\vskip1ex
\noindent
Geffroy, J., Girard, S., Jacob, P. (2006). Asymptotic normality of the $L_1-$error of a boundary estimator. {\it Nonparam. Statist.} {\bf 18}(1), 21--31.
\vskip1ex
\noindent
Girard, S., Guillou, A., Stupfler, G. (2013). Frontier estimation with kernel regression on high order moments. {\it J. Multivariate Anal.} {\bf 116}, 172--189.
\vskip1ex
\noindent
Girard, S., Iouditski, A., Nazin, A. (2005). $L_1-$optimal nonparametric frontier estimation via linear programming. {\it Autom. Remote Control} {\bf 66}(12), 2000--2018.
\vskip1ex
\noindent
Girard, S., Jacob, P. (2008). Frontier estimation via kernel regression on high power-transformed data. {\it J. Multivariate Anal.} {\bf 99}, 403--420.
\vskip1ex
\noindent
Girard, S., Jacob, P. (2009). Frontier estimation with local polynomials and high power-transformed data. {\it J. Multivariate Anal.} {\bf 100}(8), 1691--1705.
\vskip1ex
\noindent 
Hall, P. (1982). On estimating the endpoint of a distribution. {\it Ann. Statist.} {\bf 10}(2), 556--568.
\vskip1ex
\noindent
H\"ardle, W., Janssen, P., Serfling, R. (1988). Strong uniform consistency rates for estimators of conditional functionals. {\it Ann. Statist.} {\bf 16}, 1428--1449.
\vskip1ex
\noindent 
H\"ardle, W., Marron, J.S. (1985). Optimal bandwidth selection in nonparametric regression function estimation. {\it Ann. Statist.} {\bf 13}(4), 1465--1481.
\vskip1ex
\noindent
H\"ardle, W., Park, B.U., Tsybakov, A. (1995). Estimation of non-sharp support boundaries. {\em J. Multivariate Anal.} {\bf 55}, 205--218.
\vskip1ex
\noindent
Hoeffding, W. (1963). Probability inequalities for sums of bounded random variables. {\it J. Amer. Statist. Assoc.} {\bf 58}, 13--30.
\vskip1ex
\noindent
Jacob, P., Suquet, C. (1995). Estimating the edge of a Poisson process by orthogonal series. {\em J. Statist. Plann. Inference} {\bf 46}, 215--234.
\vskip1ex
\noindent 
Korostelev, A., Simar, L., Tsybakov, A. (1995). Efficient estimation of monotone boundaries. {\it Ann. Statist.} {\bf 23}, 476--489.
\vskip1ex
\noindent
Lemdani, M., Ould-Saïd, E., Poulin, N. (2009). Asymptotic properties of a conditional quantile estimator with randomly truncated data. {\it J. Multivariate Anal.} {\bf 100}, 546--559.
\vskip1ex
\noindent
Mack, Y.P., Silverman, B.W. (1982). Weak and strong uniform consistency of kernel regression estimates. {\it Z. Wahrscheinlichkeitstheorie verw. Gebiete} {\bf 61}, 405--415.
\vskip1ex
\noindent
Nadaraya, E.A. (1965). On non-parametric estimates of density functions and regression curves. {\it Theory Probab. Appl.} {\bf 10}, 186--190. 
\vskip1ex
\noindent
Parzen, E. (1962). On estimation of a probability density function and mode. {\it Ann. Math. Statist.} {\bf 33}(3), 1065--1076.
\vskip1ex
\noindent
Rosenblatt, M. (1956). Remarks on some nonparametric estimates of a density function. {\it Ann. Math. Statist.} {\bf 27}(3), 832--837.
\vskip1ex
\noindent
Silverman, B.W. (1978). Weak and strong uniform consistency of the kernel estimate of a density and its derivatives. {\it Ann. Statist.} {\bf 6}(1), 177-184.
\vskip1ex
\noindent
Stute, W. (1982). A law of the iterated logarithm for kernel density estimators. {\it Ann. Probab.} {\bf 10}, 414--422.
\end{document}